\documentclass[12pt,reqno]{amsart}

\usepackage{amsmath,amssymb,amscd,amsthm,bm}
\usepackage[all]{xy}
\usepackage[latin1]{inputenc}
\usepackage[T1]{fontenc}
\usepackage[hidelinks]{hyperref}
\usepackage{ifthen}
\usepackage{fullpage}
\usepackage{appendix}
\usepackage{url}
\usepackage{enumitem}
\usepackage{color}

\setlist{topsep=0pt,itemsep=6pt}

\newboolean{workmode}
\setboolean{workmode}{false}

\newboolean{sections}
\setboolean{sections}{true}

\newcommand{\NN}{{\mathbb{N}}}  
\newcommand{\QQ}{{\mathbb{Q}}}  
\newcommand{\RR}{{\mathbb{R}}}  
\newcommand{\Rpos}{{\mathbb{R}_{>0}}}  
\newcommand{\R}{{\mathbb{R}}}  
\newcommand{\Z}{{\mathbb{Z}}}  
\newcommand{\Q}{{\mathbb{Q}}}  
\newcommand{\ZZ}{{\mathbb{Z}}}  

\newcommand{\sgn}{{\operatorname{sgn}}}  
\newcommand{\supp}{{\operatorname{supp}}}  
\def \Cinf {C^\infty}

\newcommand{\pr}{{\operatorname{pr}}} 

\newcommand{\CIN}{{C^\infty}}   
\newcommand{\op}{{\operatorname{op}}} 
\newcommand{\Open}{{\mathbf{Open}}} 
\newcommand{\Set}{{\mathbf{Set}}} 
\newcommand{\Spec}{{\operatorname{Spec}}} 

\newcommand{\GL}{{\operatorname{GL}}}  
\renewcommand{\O}{{\operatorname{O}}} 
\newcommand{\SO}{{\operatorname{SO}}}  
\newcommand{\U}{{\operatorname{U}}}  

\DeclareMathOperator {\Lip} {Lip}

\newcommand{\ifwork}[1]{\ifthenelse{\boolean{workmode}}{#1}{}}
\newcommand{\comment}[1]{}
\newcommand{\mute}[2]{}
\newcommand{\printname}[1]{}

\ifwork{
  \renewcommand{\comment}[1]{{             \ \scriptsize{#1}\ }}
\renewcommand{\mute}[2]{{\scriptsize \ #1\ }\marginpar{\scriptsize muted}}
\renewcommand{\printname}[1]
    {{\makebox[0pt]{\hspace{-1.0in}\raisebox{8pt}{\tiny #1}}}}
}

\newcommand{\labell}[1] { \printname{#1} \label{#1} }

\newcommand{\ifsection}[2]{\ifthenelse{\boolean{sections}}{#1}{#2}}

\numberwithin{equation}{section}

\theoremstyle{plain}

\newtheorem{theorem}[equation]{Theorem}

\newtheorem{proposition}[equation]{Proposition}

\newtheorem{corollary}[equation]{Corollary}

\newtheorem{lemma}[equation]{Lemma}
\newtheorem*{lemma*}{Lemma}

\theoremstyle{definition}
\newtheorem{definition}[equation]{Definition}
\newtheorem*{definition*}{Definition}
\newtheorem{example}[equation]{Example}
\newtheorem*{example*}{Example}

\newtheorem{remark}[equation]{Remark}

\newtheorem*{remark*}{Remark}

\newtheorem{question}[equation]{Question}

\definecolor{jaw}{rgb}{0,.5,0}
\definecolor{purple}{rgb}{.6,0,.7}

\newcommand{\ynote}[1]{{\comment{\color{purple} #1}}} 
\newcommand{\jnote}[1]{{\comment{\color{jaw} #1}}} 

\def\eoe{\unskip\ \hglue0mm\hfill$\between$\smallskip\goodbreak}
\def\eor{\unskip\ \hglue0mm\hfill$\diamond$\smallskip\goodbreak}
\def\eod{\unskip\ \hglue0mm\hfill$\diamond$\smallskip\goodbreak}

\def \calO {{\mathcal O}}

\def \calI {{\mathcal I}}

\def \calC {{\mathcal C}}
\def \calF {{\mathcal F}}

\def \calD {{\mathcal D}}
\def \bfPhi {{\mathbf{\Phi}}}
\def \bfPi {{\mathbf{\Pi}}}
\def \bfXi {{\mathbf{\Xi}}}
\def \bfGamma {{\mathbf{\Gamma}}}

\def \ssminus {\smallsetminus}
\def \eps {\epsilon}
\def \del {\partial}
\def \ol  {\overline}
\def \ul  {\underline}

\makeatletter
\def\th@plain{%
  \thm@notefont{}
  \itshape 
}
\def\th@definition{%
  \thm@notefont{}
  \normalfont 
}
\makeatother



\begin{document}

\title{Diffeological, Fr\"olicher, and Differential Spaces}
\date{\today}

\author{Augustin Batubenge}
\email{a.batubenge@gmail.com}

\author{Yael Karshon}
\address{Department of Mathematics, University of Toronto, Toronto, ON, Canada; and School of Mathematical Sciences, Tel-Aviv University, Tel-Aviv, Israel}
\email{karshon@math.toronto.edu}
\email{yaelkarshon@tauex.tau.ac.il}

\author{Jordan Watts}
\address{Department of Mathematics, Central Michigan University, Pearce Hall 214, Mount Pleasant, MI 48859}
\email{jordan.watts@cmich.edu}

\begin{abstract}
Differential calculus on cartesian spaces has many generalisations.  
In particular,
on a set $X$, a diffeological structure is given 
by maps from open subsets of cartesian spaces to $X$,
a differential structure is given by maps from $X$ to $\RR$, and
a Fr\"olicher structure is given by 
maps from $\RR$ to $X$ as well as maps from $X$ to $\RR$.
We illustrate the relations between these structures through examples.
\end{abstract}

\maketitle

\section{Introduction}
\labell{sec:intro}

There are many structures in the mathematical literature 
that generalise differential calculus beyond manifolds.
In this paper we focus on the simplest such structures:
diffeology (as defined by Souriau), 
differential structures (in the sense of Sikorski), 
and Fr\"olicher structures.
A diffeology on a set $X$ is given by a set of maps from open subsets
of cartesian spaces to $X$; see Definition~\ref{d:diffeology}.
A differential structure on a set $X$
is given by a set of maps from $X$ to $\RR$;
see Definition~\ref{d:diff space}.
A Fr\"olicher structure on a set $X$
is given by a set of maps from $\RR$ to $X$
and a set of maps from $X$ to $\RR$;
see Definition~\ref{d:frolicher}.
These structures are motivated by the following characterisations
of smooth maps between manifolds.

Let $M$ and $N$ be open subsets of cartesian spaces $\RR^m$ and $\RR^n$
and $\psi \colon M \to N$ a function.
Smoothness of $\psi$ is equivalent to each of the following conditions.
\begin{enumerate}[topsep=0pt]
\item
For each $k$,
each open subset $U$ of $\RR^k$,
and each smooth map $p \colon U \to M$,
the composition $\psi \circ p \colon U \to N$ is smooth.
\item
For each real-valued smooth function $f \colon N \to \RR$,
the composition $f \circ \psi \colon M \to \RR$ is smooth.
\item
For each smooth curve $\gamma \colon \RR \to M$,
the composition $\psi \circ \gamma \colon \RR \to N$ is smooth.
\end{enumerate}
The fact that the third condition implies the smoothness of $\psi$
follows from the following theorem of Jan Boman \cite[Theorem 1]{boman}:
Let $f$ be a function from $\RR^d$ to $\RR$,
and assume that the composition $f \circ u$ is in $\CIN(\RR,\RR)$
for every $u \in \CIN(\RR,\RR^d)$. Then $f$ is in $\CIN(\RR^d,\RR)$.

In this paper, we illustrate the
relation between differential structures, diffeological structures,
and Fr\"olicher structures, through examples.  
The paper should be valuable
to researchers and graduate students seeking a quick and effective 
introduction to these structures.  
Whereas the goal of the paper is mostly expository,
the paper does contain new material.
We identify Fr\"olicher spaces with so-called reflexive differential spaces 
and so-called reflexive diffeological spaces
(see Definition~\ref{d:reflexive} and Theorems \ref{t:A} and \ref{t:B}).
This notion of reflexivity
that we introduce and its relation to Fr\"olicher spaces
was in principle known to experts but to our knowledge
it has not been made explicit in the literature.
Throughout the paper
we place particular emphasis on the theme of reflexivity and non-reflexivity.
Our examples --- many of which are new ---
illustrate the vast richness of diffeology and differential
structures beyond the reflexive ones,
while pointing at some that are reflexive for non-trivial reasons.
Last but not least, we pose a number of open questions.

We deliberately focus on these structures, 
which we view as the simplest among the many
generalisations of differential calculus.  To this end, we do not focus
on higher categorical approaches to smoothness such as differentiable
stacks, nor algebro-geometric settings such as $\CIN$-schemes, nor
differentiability of finite order. We believe that a good understanding
of the simpler structures would be beneficial also for those who
wish to work with other generalisations of differential calculus,
as different tools capture different subsets of the phenomena that
we illustrate.  
Nevertheless, in response to a referee request,
we are including an appendix (Appendix~\ref{app:comparisons})
in which we comment on relations of these simpler structures
to some other structures in the higher categorical and algebro-geometric settings.

The richness of these non-reflexive examples
motivates working in the presence of both a diffeology and a 
differential structure that are compatible but not necessarily reflexive.
Such spaces, named ``Watts spaces'' in \cite{virgin},
have been informally promoted by Jordan Watts for many years.

In Section~\ref{s:maintheorems},
we identify the category of Fr\"olicher spaces
with the categories of so-called reflexive differential spaces
and so-called reflexive diffeological spaces
(see Definition~\ref{d:reflexive}
and Theorems \ref{t:A} and \ref{t:B}) 
and give some examples of non-reflexive diffeological spaces 
and non-reflexive differential spaces.

One of the motivations for considering diffeological and differential
structures is that they are meaningful for arbitrary subsets and
quotients of manifolds.  In Section~\ref{sec:subsets}, we discuss how
diffeological and differential structures relate on these objects; see
Propositions~\ref{p:subsets} and~\ref{p:quotients}.  
We include two open questions, one about the diffeology
of symplectic reduced spaces, and one about the differential structure
of an irrational line in the torus.

In Section~\ref{sec:manifolds}, we consider orbifolds, quotients by
compact Lie group actions, and manifolds with corners.  By a result
of Gerald Schwarz \cite{schwarz}, the \emph{Hilbert map} identifies the
quotient of a linear compact Lie group action with a subset of a cartesian
space as differential spaces.  Consequently, 
the subspace differential structure 
on its image is reflexive.  
As a consequence, manifolds with corners can be defined
equivalently as differential spaces or as diffeological spaces; either
of these structures is reflexive.  See Example~\ref{x:corners}.  On the
other hand, the quotient diffeology can be non-reflexive, and consequently
different from the subset diffeology on the image of the Hilbert map.
See Examples~\ref{x:orthogonal action} and \ref{x:orbifolds}.

In Section~\ref{sec:wedge}, we consider finite unions of copies of
the real line, as well as finite unions of manifolds.
For example, the union of the three coordinate axes in
$\RR^3$ is diffeomorphic to a union of three concurrent lines in $\RR^2$
diffeologically but not as differential spaces.  
The former differential space is reflexive; the latter is not.
See Examples~\ref{x:wedgesum} and \ref{three lines}.
For a generalization to cleanly intersecting submanifolds,
see Example~\ref{x:wedgesum2}.
We conclude this section with a question
about the reflexivity of an example that comes from polyfolds.

In Section~\ref{sec:topology}, we consider topological properties
of diffeological and differential spaces, and obtain some 
topological necessary conditions for reflexivity.

Appendix~\ref{proofs} contains some technical proofs
that are deferred from the earlier sections.

In Appendix~\ref{app:comparisons},
we briefly compare diffeological, differential, and Fr\"olicher
spaces to Lie groupoids, stacks, sheaves of sets over the site $\Open$, Mostow spaces, subcartesian spaces, differentiable
spaces \`a la Gonzalez-Salas/Spallek, and $\CIN$-schemes.

This paper evolved from visits of Patrick Iglesias-Zemmour and author Batubenge
to the University of Toronto.
Batubenge and Iglesias-Zemmour contributed to the initiation
and vision of this paper and to Sections~\ref{s:maintheorems}, 
\ref{sec:subsets}, \ref{sec:manifolds}, \ref{sec:wedge}, 
and relevant proofs in Appendix~\ref{proofs}.
Many of the details were worked out and written up by Watts
as Chapter 2 of his University of Toronto Ph.D.\ thesis \cite{watts},
supervised by Yael Karshon.
Appendix~\ref{app:comparisons} was authored mostly by Jordan Watts,
in response to our referee's request to clarify how the simple structures 
on which we are focusing in this paper
relate to more complicated structures that occur in the literature.

\smallskip\noindent\textbf{Some history and notes on the literature.}

The development of the various notions of smooth structures discussed in this paper occurred mainly in the 1960s, '70s, and '80s, motivated by the need to push differentiability beyond the confines of finite-dimensional manifolds to the singular subset, singular quotient, and infinite-dimensional settings.

Differential structures (Definition~\ref{d:diff space})
were introduced by Sikorski in the late 1960s;
see \cite{sikorski1,sikorski2}.  Many of the properties of the smooth structure on a smooth manifold can be derived from its ring of smooth functions; a differential space is a topological space equipped with a ring of functions that captures these properties.  A differential structure determines a sheaf of continuous functions that contains the constants (as considered by Hochschild \cite{hochschild}), which, in turn, is a special case of a ringed space (a topological space equipped with a sheaf of rings; see EGA 1 \cite{grothendieck}).  Pushing similar notions from algebraic geometry into the realm of differential geometry leads to further developments, $\CIN$-schemes \cite{joyce} and differentiable spaces in the sense of Gonzalez-Salas \cite{GS} being some resulting theories.

Special cases of differential spaces appear in the literature in various contexts.  A subcartesian space, introduced by Aronszajn in the late 1960s and motivated by manifolds with
singularities that occur in his study of the Bessel potential in functional analysis, is a Hausdorff differential space that is locally diffeomorphic
to (arbitrary) subsets of cartesian spaces; see \cite{A,AS1,AS2,sniatycki}.  In the mid-1970s, interest in equipping singular orbit spaces of compact Lie group actions with a smooth structure (see Bredon \cite[Chapter 6]{bredon}) led to a result of Schwarz \cite{schwarz} showing that while \emph{a priori} quotient spaces, these spaces are in fact subcartesian as well; see also Cushman-\'Sniatycki \cite{CS}, who work with orbit spaces of proper Lie group actions in the subcartesian setting.  A similar result for symplectic quotients in the early 1990s by Arms-Cushman-Gotay \cite{ACG} lead to the study of these spaces as subcartesian spaces equipped with Poisson structures; for example, this is used in the treatment of symplectic quotients as symplectic stratified spaces by Sjamaar-Lerman \cite{sjamaar-lerman}.  Today, viewing stratified spaces as differential spaces is commonplace; for example, they appear in the book by Pflaum \cite{pflaum} and the work of \'Sniatycki \cite{sniatycki} as subcartesian spaces, and Kreck's stratifolds \cite{kreck,tene} are a version of stratified differential spaces whose functions satisfy certain conditions at the strata.  Differential subspaces of cartesian spaces are special cases of subcartesian spaces: convex subsets are studied in Karshon-Watts \cite{KW:convex}, and so-called H\"older sets and certain subanalytic sets are studied in Rainer \cite[Theorem 1.13]{rainer}.
Our main reference on the theory of differential spaces is the book by \'Sniatycki \cite{sniatycki}.

Diffeology (Definition~\ref{d:diffeology}) was introduced
by Jean-Marie Souriau around 1980; see \cite{souriau}.  
An early success of the theory which helped to motivate its further
development is the work of Donato-Iglesias on irrational tori
\cite{donato-iglesias}; see Example~\ref{x:irrational flow}.
Irrational tori, (or ``infracircle''s) 
also appear in the study of geometric quantisation 
\cite{weinstein}, \cite{iglesias:95}, as well as the integration 
of certain Lie algebroids \cite{crainic}, and so play a role 
in mathematical physics.
Stratified spaces as diffeological spaces appear in 
\cite{GIZ1,GIZ2},
where the diffeological language is a natural setting to describe the
so-called zero perverse differential forms of the intersection theory
of Goresky-MacPherson. 
Our main reference on the theory of diffeology
is the book by Iglesias-Zemmour \cite{iglesias}.

Souriau's motivation for developing diffeology came from
infinite-dimensional groups appearing in mathematical physics.
A similar notion was introduced and studied by Kuo-Tsai Chen already
in the 1970s for the purpose of putting differentiability on path
spaces used in variational calculus on an equal footing with smooth
structures on manifolds; the precise definition went through several
revisions \cite{chen1,chen2,chen3,chen4}.  The main difference
between diffeological spaces and Chen spaces is that the latter use
convex subsets instead of open subsets of cartesian spaces as domains
of the so-called plots. In \cite{KW:convex}, authors Karshon and Watts show that diffeological spaces are isomorphic as a category to a full subcategory of Chen spaces.

Similar motivations in functional analysis lead Fr\"olicher and Kriegl 
to introduce what are now called Fr\"olicher spaces in their book \cite{FK}, 
following the work of Fr\"olicher in the early 1980s 
\cite{frolicher2,frolicher,frolicher3}.  A special case is the
``convenient setup'' of Fr\"olicher, Kriegl, and Michor \cite{FK,KM},
which applies to finite and infinite-dimensional vector spaces and
manifolds.  Vector spaces from the diffeological perspective are studied by Christensen-Wu in \cite{CW:vector spaces}.  Reflexivity of spaces of smooth maps between diffeological spaces is examined in an upcoming paper \cite{karshon-watts:part2} by the authors Karshon and Watts.
Iglesias-Zemmour and Karshon study Lie groups as diffeological subgroups of diffeomorphism groups in \cite{IZK}, coadjoint orbits of infinite-dimensional groups are studied by Iglesias-Zemmour in \cite{donato-iglesias:coadjoint} and Lee in \cite{BrianLee}, infinite products and coproducts appear in Karshon's paper on moduli spaces \cite{karshon:moduli}, diffeological classifying spaces appear in the work of Magnot-Watts \cite{magnot-watts} and Christensen-Wu \cite{CW:classifying}, and Magnot studies Fr\"olicher and diffeological Lie groups in \cite{magnot:ambrose-singer, magnot:diffeos, magnot-reyes}.

Many of the categories mentioned above are compared in the paper 
of Andrew Stacey \cite{stacey}.  Along with the diffeological,
differential, and Fr\"olicher spaces, he also considers various
definitions of Chen spaces, as well as Smith spaces \cite{smith}
(topological spaces equipped with a set of continuous functions that
satisfy a certain ``reflexivity'' condition), and constructs functors 
between these categories.  Treatments of diffeological and
Chen spaces from the point-of-view of sheaves on categories is given
in Baez-Hoffnung \cite{BH} (also see \cite{cafe}), and a treatment of
diffeological spaces from the point-of-view of stacks on manifolds is
given in Watts-Wolbert \cite{WW}.

Quotient spaces of Lie group actions, Lie groupoids, and more generally, singular foliations, form another setting in which the theories of diffeology, Fr\"olicher spaces, and differential spaces are important.  Orbifolds are given a diffeological treatment in a paper of Iglesias-Zemmour, Karshon, and Zadka \cite{IZKZ}, and further studied in an intersection of diffeology with non-commutative geometry in a paper of Iglesias-Zemmour and Laffineur \cite{iglesias-laffineur}.  Quasifolds from a diffeological and groupoid perspective are treated in Karshon-Miyamoto \cite{quasifolds}, based on the work of Masrour Zoghi \cite{zoghi}.  From a differential space perspective, or equivalently in this case, a Fr\"olicher perspective, orbifolds are examined in a paper of Watts \cite{watts-orb}, and orbit spaces of linear circle actions in a paper of Craig-Downey-Goad-Mahoney-Watts \cite{CDGMW}.  Differential forms on these quotient spaces from a diffeological perspective are compared with basic forms for compact Lie group actions in Watts' Ph.D. thesis \cite{watts}, proper Lie group actions in a paper by Karshon-Watts \cite{karshon-watts:basic}, proper Lie groupoids in a paper by Watts \cite{watts-gpds}, and certain singular foliations by Miyamoto \cite{miyamoto}.

\subsection*{Acknowledgements}
The authors would like to thank Patrick Iglesias-Zemmour for his many contributions to this project.
Yael Karshon would like to thank Katrin Wehrheim,
Shintaro Kuroki, Eugene Lerman, 
and Haggai Tene for helpful discussions, 
and Peter Michor for a helpful discussion back in 2003.
Yael Karshon and Jordan Watts are grateful to Enxin Wu 
for helpful discussions.
Jordan Watts would like to thank Jo\~ao Nuno Mestre for a helpful discussion
regarding the contents of Appendix~\ref{app:comparisons}.

This research is partially supported by the Natural Sciences and
Engineering Research Council of Canada.
Augustin Batubenge would also like to thank CSET-CRIC,
Topology at Unisa Flagship, and Professor Fran\c{c}ois Lalonde
for their support.

\section{Relations Between Structures}
\labell{s:maintheorems}

\begin{definition}[Diffeology] \ 
\labell{d:diffeology}
Let $X$ be a nonempty set.  A \emph{parametrisation} of $X$ is a function
$p \colon U \to X$ where $U$ is an open subset of $\RR^n$
for some $n$.  A \emph{diffeology} $\mathcal{D}$ on $X$ is a set of
parametrisations satisfying the following three conditions.
\begin{enumerate}
\item \textbf{(Covering)}
For every $x\in X$ and every nonnegative integer
$n\in\NN$, the constant function $p \colon \RR^n\to\{x\}\subseteq X$ is in
$\mathcal{D}$.
\item \textbf{(Locality)}
Let $p \colon U \to X$ be a parametrisation such that for
every $u\in U$ there exists an open neighbourhood $V$ of $u$ in $U$
satisfying $p|_V\in\mathcal{D}$. Then $p\in\mathcal{D}$.
\item \textbf{(Smooth Compatibility)}
Let $p \colon U\to X$ be a plot in $\mathcal{D}$.
Then for every $n\in\NN$, every open subset $V\subseteq\RR^n$, and every
infinitely-differentiable map $F \colon V\to U$, 
we have $p\circ F\in\mathcal{D}$.
\end{enumerate}
A set $X$ equipped with a diffeology $\mathcal{D}$ is called a
\emph{diffeological space} and is denoted by $(X,\mathcal{D})$.
When the diffeology is understood, we may drop the symbol $\calD$.
The parametrisations in $\mathcal{D}$ are called \emph{plots}.
A map $F \colon X \to Y$ between diffeological spaces
is \emph{smooth} 
if for any plot $p \colon U \to X$ of $X$
the composition $F \circ p \colon U \to Y$ is a plot of $Y$.
The map is a \emph{diffeomorphism} if it is smooth and has a smooth inverse.
To avoid ambiguity, we sometimes say that the map is
\emph{diffeologically smooth} or is a \emph{diffeological diffeomorphism}.

The \emph{D-topology} on $X$ is the strongest topology
in which every plot is continuous;
thus, a subset $Y$ of $X$ is D-open if and only if for each plot $p \in \calD$
the preimage $p^{-1}(Y)$ is open in the domain of $p$.
\eod
\end{definition}

Given a collection of functions $\calF_0$ on a set $X$,
its \emph{initial topology} is the weakest topology on $X$ 
for which every function in $\calF_0$ is continuous.
Thus, a sub-basis for the initial topology 
is given by the pre-images of open intervals by functions in $\calF_0$.

\begin{definition}[Differential space] \labell{d:diff space}
Let $X$ be a nonempty set.  A \emph{differential structure} on $X$ 
is a nonempty family $\mathcal{F}$ of real-valued functions on $X$,
along with its initial topology, satisfying the following two conditions.
\begin{enumerate}
\item \textbf{(Smooth compatibility)} For any positive integer $k$, functions
$f_1,...,f_k\in\mathcal{F}$, and $F\in\CIN(\RR^k)$, the composition
$F(f_1,...,f_k)$ is in $\mathcal{F}$.
\item \textbf{(Locality)} Let $f \colon X\to\RR$ be a function such that
for any $x\in X$
there exist an open neighbourhood $U\subseteq X$ of $x$ and a function
$g\in\mathcal{F}$ satisfying $f|_U=g|_U$.  Then $f\in\mathcal{F}$.
\end{enumerate}
A set $X$ equipped with a differential structure $\mathcal{F}$ is called
a \emph{differential space} and is denoted by $(X,\mathcal{F})$.
When the differential structure is understood, we may drop the symbol $\calF$.
A map $F \colon X\to Y$ between differential spaces $(X,\mathcal{F}_X)$
and $(Y,\mathcal{F}_Y)$ is \emph{smooth}
if for every function $f \colon Y \to \R$ in $\calF_Y$
the composition $f \circ F$ is in $\calF_X$.
The map is a \emph{diffeomorphism}
if it is smooth and has a smooth inverse.
To avoid ambiguity, we sometimes say that the map is
\emph{functionally smooth} or is a \emph{functional diffeomorphism}.
\eod
\end{definition}

\begin{definition}[``Compatible'' and ``induces'']
\labell{d:compatible induce}
Given a set $X$ with a collection $\calD_0$ of parametrisations
and a collection $\calF_0$ of real-valued functions, we say that 
\begin{enumerate}
\item[(i)] 
$\calD_0$ and $\calF_0$ are \emph{compatible} if $f \circ p$ is
infinitely-differentiable
for all $p \in \calD_0$ and $f \in \calF_0$;
\item[(ii)]
$\calD_0$ \emph{induces} $\calF_0$ if $\calF_0$ coincides 
with the set
$$   \Phi\mathcal{D}_0 := \{ f \colon X \to \RR \ | \
   \forall (p \colon U \to X) \in \mathcal{D}_0,
   \ f\circ p \in \CIN(U) \} $$
of those real-valued functions 
whose precomposition with each element of $\calD_0$ 
is infinitely-differentiable;
\item[(iii)]
$\calF_0$ \emph{induces} $\calD_0$
if $\calD_0$ coincides with the set 
$$
   \Pi \mathcal{F}_0 := \{
   \text{parametrisations } p \colon U \to X \ | \
   \forall f \in \mathcal{F}_0, \ f \circ p \in \CIN(U) \} $$
of those parametrisations
whose composition with each element of $\calF_0$ is infinitely-differentiable.
\end{enumerate}
Thus, $\calD_0$ and $\calF_0$ are compatible
if and only if $\calF_0$ is contained in $\Phi\calD_0$,
if and only if $\calD_0$ is contained in $\Pi\calF_0$.
\eod
\end{definition}

\begin{example}[Manifolds] \labell{e:manifolds}
On a smooth manifold $M$, the sets of parametrisations $U \to M$ 
that are infinitely-differentiable and the set of real-valued functions 
$M\to\RR$ that are infinitely-differentiable are a diffeology and a
differential structure that induce each other.  This follows from
the fact that smoothness is a local property
and from the existence of smooth bump functions.
\eoe
\end{example}

\begin{remark} \labell{properties}
We make the following easy observations:
\begin{itemize}
\item
Each of the operations $\calD_0 \mapsto \Phi \calD_0$
and $\calF_0 \mapsto \Pi \calF_0$ is inclusion-reversing.
\item
We always have $\Pi \Phi \calD_0 \supseteq \calD_0$
and $\Phi \Pi \calF_0 \supseteq \calF_0$.
\end{itemize}
These facts imply that, given a family $\calD$ of parametrisations,
there exists a family of real-valued functions that induces $\calD$ 
if and only if $\Pi\Phi\calD=\calD$.
Indeed, if $\calD=\Pi\calF$ then $\Pi\Phi\calD \subseteq \calD$
amounts to $\Pi\Phi\Pi\calF \subseteq \Pi\calF$,
which follows from $\Phi\Pi\calF \supset \calF$.
Similarly, given a family $\calF$ of real-valued functions,
there exists a family of parametrisations 
that induces $\calF$ if and only if $\Phi\Pi\calF=\calF$.
\eor
\end{remark}

\begin{definition}[Reflexive]
\labell{d:reflexive}
A diffeology $\mathcal{D}$ is \emph{reflexive} if 
$\Pi\Phi\mathcal{D}=\mathcal{D}$.  
A differential structure $\mathcal{F}$ is \emph{reflexive} if
$\Phi\Pi\mathcal{F}=\mathcal{F}$.
\eod
\end{definition}

\begin{proposition}[Reflexive stability] \labell{p:reflexive}
For any family $\calF_0$ of real-valued functions on a set,
\ $\Pi\calF_0$ is a reflexive diffeology on the set.
For any family $\calD_0$ of parametrisations on a set,
\ $\Phi\calD_0$ is a reflexive differential structure on the set.
\end{proposition}
We prove Proposition \ref{p:reflexive} in \S\ref{subsec:reflexive stability}.

Thus, if a diffeology $\calD$ and a differential structure $\calF$
induce each other, then they are both reflexive.
For example, manifolds are reflexive both as diffeological spaces
and as differential spaces.
Here are examples of diffeological and of differential spaces
that are not reflexive:

\begin{example}[Spaghetti diffeology]
\labell{wire}
The \emph{spaghetti diffeology} (or wire diffeology) on $\RR^2$
consists of those parametrisations that locally factor through curves.
That is, a parametrisation $p \colon U \to \RR^2$ is in the spaghetti diffeology
if and only if for every $u\in U$ there is
an open neighbourhood $V$ of $u$ in $U$,
a smooth map $F \colon V \to \RR$,
and a smooth curve $q \colon \RR \to \RR^2$, such that
$ p|_V = q \circ F$.
See \cite[Section 1.10]{iglesias} for more details.

The differential structure that is induced by the spaghetti diffeology
consists of those real-valued functions $f \colon \RR^2 \to \RR$
such that $f \circ q$ is smooth for every smooth curve $\RR \to \RR^2$.
By Boman's theorem \cite[Theorem 1]{boman},
every such function $f$ is infinitely-differentiable.
Thus, this is the standard differential structure on $\RR^2$,
and the diffeology that it induces is the standard diffeology on $\RR^2$.

The spaghetti diffeology and the standard diffeology have the same smooth curves
$\R \to \R^2$, but they are different.
For example, the identity map on $\RR^2$ is in the standard diffeology 
but not in the spaghetti diffeology.  Thus, the spaghetti diffeology is not reflexive.
\eoe
\end{example}

\begin{example}[Rational numbers] \labell{e:rational numbers}
Consider the set $\QQ$ of rational numbers 
with the differential structure $\CIN(\QQ)$
that consists of those functions $f \colon \QQ \to \RR$
that locally extend to smooth functions on $\RR$.
This includes, for example, the restriction to $\QQ$ of the function
$x \mapsto \frac{1}{x-\sqrt{2}}$.  All the plots in $\Pi\CIN(\QQ)$ are locally constant.
(Indeed, since the inclusion map $\QQ\hookrightarrow \RR$ is in $\CIN(\QQ)$, every $p\in\Pi\CIN(\QQ)$ must be smooth as a function to $\RR$.  By the intermediate value theorem, such a $p$ must be locally constant.)  Consequently, the differential space $(\QQ,\CIN(\QQ))$ is not reflexive.
\eoe
\end{example}

\begin{example}[$\mathbf{C^k(\RR)}$]\labell{Ck}
Fix an integer $k\geq 0$.  Consider the real line $\RR$ with the differential structure $C^k(\RR)$ consisting of those real-valued functions that are $k$-times continuously differentiable.  All the plots in $\Pi C^k(\RR)$ are locally constant. (Indeed, take any parametrisation $p\colon U\to\RR$.  Since the identity map is in $C^k(\RR)$, if $p\in\Pi C^k(\RR)$, then $p$ must be infinitely-differentiable.  If $p$ is infinitely-differentiable and not locally constant, then there exists $u\in U$ such that $dp|_u\neq 0$; the composition of $p$ with a map $f\in C^k(\RR)$ that is not smooth at $p(u)$ is not smooth, so $p\notin\Pi C^k(\RR)$.)  Consequently, the differential space $(\RR,C^k(\RR))$ is not reflexive.
\eoe
\end{example}

Diffeological spaces, along with diffeologically smooth maps, form a category;
reflexive diffeological spaces form a full subcategory.
Differential spaces, along with functionally smooth maps, form a category;
reflexive differential spaces form a full subcategory.

If $(X,\mathcal{D}_X)$ and $(Y,\mathcal{D}_Y)$ are two diffeological spaces
and $F \colon X \to Y$ is a diffeologically smooth map,
then $F$ is also a functionally smooth map
from $(X,\Phi\mathcal{D}_X)$ to $(Y,\Phi\mathcal{D}_Y)$.
Thus, we have
a functor $\mathbf{\Phi}$ from diffeological
spaces to reflexive differential spaces that sends
a diffeological space $(X,\mathcal{D})$
to the reflexive differential space $(X,\Phi\mathcal{D})$
and that sends each map to itself.  Similarly, we have a functor $\mathbf{\Pi}$
from differential spaces to reflexive
diffeological spaces that sends a differential space $(X,\mathcal{F})$
to the reflexive diffeological space $(X,\Pi\mathcal{F})$
and that sends each map to itself.
In \S\ref{subsec:proof of A} we prove these facts
and obtain the following theorem:

\begin{theorem}
[Isomorphism of categories of reflexive spaces]
\labell{t:A}
The restriction of the functor $\bfPhi$ to the subcategory
of reflexive diffeological spaces is an isomorphism of categories
onto the subcategory of reflexive differential spaces.  The restriction of the functor $\bfPi$ to the subcategory
of reflexive differential spaces is an isomorphism of categories
onto the subcategory of reflexive diffeological spaces.
These isomorphisms are inverses of each other.
\end{theorem}

Given a set $X$ and a family $\mathcal{F}_0$ of real-valued functions on $X$,
we also consider the set $\Gamma \calF_0$ of those maps from $\R$ to $X$
whose composition with each element of $\calF_0$ is infinitely-differentiable: 
$$
 \Gamma\mathcal{F}_0:= \{ c \colon \RR\to X~|~\forall
 f\in\mathcal{F}_0,~f\circ c\in\CIN(\RR)\}.
$$
The operation $\calF_0 \mapsto \Gamma\calF_0$ is inclusion-reversing.
Also, for any family of functions $\calF_0$ from $X$ to $\RR$
and family of functions $\calC_0$ from $\RR$ to $X$, we have
$\calC_0 \subseteq \Gamma\Phi\calC_0$ and $\calF_0 \subseteq \Phi\Gamma\calF_0$.
These facts imply that $\Gamma \Phi \Gamma \calF_0 = \Gamma \calF_0.$

\begin{definition}[Fr\"olicher spaces]\labell{d:frolicher}
A \emph{Fr\"olicher structure} on a set $X$ is a family $\mathcal{F}$ 
of real-valued functions $X \to \RR$
and a family $\mathcal{C}$ of maps $\RR\to X$,
such that
$$ \Phi\mathcal{C}=\mathcal{F} \quad \text{and} \quad 
   \Gamma\mathcal{F}=\mathcal{C}.$$  
Such a triple
$(X,\mathcal{C},\mathcal{F})$ is a \emph{Fr\"olicher space}.

Let $(X, \mathcal{C}_X, \mathcal{F}_X)$ and
    $(Y, \mathcal{C}_Y, \mathcal{F}_Y)$ be Fr\"olicher spaces.
A map $F \colon X \to Y$ 
is \emph{Fr\"olicher smooth} if it satisfies one, hence all,
of the following equivalent conditions:
\begin{itemize}
\item[(i)]
$f \circ F \in \mathcal{F}_X$ for every $f \in \mathcal{F}_Y$.
\item[(ii)]
$f \circ F \circ c \in \CIN(\RR,\RR)$
for every $c \in \calC_X$ and $f \in \calF_Y$.
\item[(iii)]
$F \circ c \in \mathcal{C}_Y$ for every $c\in\mathcal{C}_X$.
\end{itemize}
((i) implies (ii) because $\calF_X$ and $\calC_X$ are compatible,
(ii) implies (i) because $\calF_X = \Phi\calC_X$,
(ii) implies (iii) because $\calC_Y = \Gamma\calF_Y$,
(iii) implies (ii) because $\calF_Y$ and $\calC_Y$ are compatible.)
\eod
\end{definition}

Fr\"olicher spaces, along with Fr\"olicher smooth maps, form a category.
There is a functor $\bfXi$ 
from the category of Fr\"olicher spaces
to the category of reflexive differential spaces 
that takes $(X,\calC,\calF)$ to $(X,\calF)$
and takes each map to itself.
There is also a functor $\bfGamma$ 
from the category of differential spaces
that takes $(X,\calF)$ to $(X,\Gamma\calF,\Phi\Gamma\calF)$
and takes each map to itself.
In \S\ref{subsec:proof of B} we prove these facts
and obtain the following theorem:

\begin{theorem} 
[Fr\"olicher spaces as reflexive spaces]
\labell{t:B}
The functor $\bfXi$ is an isomorphism
from the category of Fr\"olicher spaces
to the category of reflexive differential spaces.
The functor $\bfGamma$ restricts to an isomorphism
from the category of reflexive differential spaces
to the category of Fr\"olicher spaces.
These isomorphisms are inverses of each other.
\end{theorem}

To summarise, we have isomorphisms between the categories of 
Fr\"olicher spaces $\{(X,\calC,\calF)\}$, 
reflexive differential spaces $\{(X,\calF)\}$,
and reflexive diffeological spaces $\{(X,\calD)\}$,
where the functors send every map to itself
and their actions on objects are given by the following
commuting diagram.
\begin{equation} \labell{six isomorphisms}
\xymatrix{ \left\{ (X,\calD) \right\} 
       \ar@/^/[rr]^{\calC=\text{1-dim'l plots},\ \calF=\Phi\calD}
       \ar@/^/[rdd]^{\calF=\Phi\calD}
 & & \{(X,\calC,\calF)\} \ar@/^/[ll]^{\calD=\Pi\calF}
                   \ar@/^/[ldd]^{\text{ same } \calF} \\ && \\
 & \{ (X,\calF) \}
                   \ar@/^/[ruu]^{\calC=\Gamma\calF}
                   \ar@/^/[luu]^{\calD=\Pi\calF} &
}
\end{equation}

\smallskip\noindent\textbf{Notes.}

\begin{enumerate}

\item
In the literature, what we call \emph{differential structure},
\emph{differential space}, \emph{functionally smooth map},
and \emph{functional diffeomorphism},
are sometimes called \emph{Sikorski structure},
\emph{Sikorski space}, \emph{Sikorski smooth map},
and \emph{Sikorski diffeomorphism}.

\item  
In the literature, the adjective ``reflexive'' 
often refers to a Banach space $E$
and means that the natural inclusion of $E$ into $(E^*)^*$ is an isomorphism.  
Many Banach spaces (for example $C([0,1])$) are not reflexive
as Banach spaces,
but the diffeology and differential structure on a Banach 
(or Fr\'echet) space
that consist of those parametrisations and those real-valued functions 
that are smooth in the usual sense
are always reflexive; see \cite{frolicher2,hain};
also see \cite{karshon-watts:part2}.
Also, the analogue of reflexive stability (Proposition~\ref{p:reflexive}) for the functor sending a Banach space to its dual is not true:  by the Hahn-Banach theorem, a Banach space $E$ is reflexive if and only if its dual space $E^*$ is reflexive \cite{folland}.

\item
The behaviour of the functors $\Phi$ and $\Pi$ is that of an \emph{antitone Galois connection} \cite{ore}.  Other examples of such relationships include sets of polynomials and their zero sets in algebraic geometry, as well as field extensions and their Galois groups.

\item
The functor $\bfXi \colon (X,\calC,\calF) \mapsto (X,\calF)$
from Fr\"olicher spaces to differential spaces
was described in Cherenack's paper \cite{CF}.
The functor $\bfGamma \colon (X,\calF) \mapsto (X,\Gamma\calF,\Phi\Gamma\calF)$
from differential spaces to Fr\"olicher spaces
was described in Batubenge's Ph.D.\ thesis \cite[\S 2.7]{batubenge-thesis}.  A differential space $(X,\calF)$ is reflexive if and only if $\Phi\Gamma\calF = \calF$; these spaces were introduced in \cite[\S 5.2]{batubenge-thesis} under the name ``pre-Fr\"olicher spaces''.  Further comparisons between Fr\"olicher and differential spaces appear in \cite{BT}.

\item Example~\ref{Ck} appears in \cite[Example 2.79]{watts}.
In the context of Smith spaces, $(\RR,C^0(\RR))$ 
 is discussed in \cite[p.100, paragraph on ``Smith spaces'']{stacey}; 
however, when $\RR$ is equipped with its standard topology,
$(\RR,C^0(\RR))$ is not a Smith space.

\item
Some of the results of this section can be rephrased in terms 
of adjoint functors and reflective subcategories; 
see Sections 8.4.1 and 8.4.4 of \cite{FK}.  
In particular, $\Xi$ is a left adjoint to $\Gamma$, $\Phi$ is a left adjoint to $\Pi$, and $\Gamma\circ\Phi$ is a left adjoint to $\Pi\circ\Xi$.
These facts are also in Stacey's paper \cite{stacey}, 
noting that $(X,\calC,\calF)$ should be $(X,\calC_X,\calF_X)$ 
in the last sentence of the second paragraph of the subsection 
on Smith and Fr\"olicher spaces (Section 5).

\end{enumerate}

\section{Subsets and Quotients}
\labell{sec:subsets}

In this section we discuss subsets and quotients
from several points of view:
differential structures, diffeology, and topology.
We omit many of the proofs.
The interested reader can fill in the details as an exercise or look them up in 
Iglesias--Zemmour's book \cite[Chapter 1]{iglesias}, 
Watts' thesis \cite[Chapter 2]{watts}, 
or \'Sniatycki's book~\cite[Chapter 2]{sniatycki}.

\begin{definition}[Subsets]\labell{d:subsets}
Let $X$ be a set and $Y \subseteq X$ a subset.
Given a diffeology $\calD$ on $X$,
the \emph{subset diffeology} on $Y$
consists of those parametrisations $p \colon U \to Y$
whose composition with the inclusion map $Y \hookrightarrow X$
is a plot in $\calD$.
Given a differential structure $\calF$ on $X$,
the \emph{subspace differential structure} on $Y$
consists of those functions $f \colon Y \to \RR$
that locally extend to $X$ in the following sense:
for every $x \in Y$ there exists an open neighbourhood $U$
of $x$ in $X$ with respect to the initial topology 
and a function $\tilde{f} \in \calF$ such that
$f|_{U \cap Y} = \tilde{f}|_{U \cap Y}$.
\eod
\end{definition}

Differential structures are well adapted to subsets:

\begin{proposition}[Differential subspaces]\labell{p:subsets}
Given a differential space $(X,\calF)$ and a subset $Y \subset X$,
we obtain on $Y$ an unambiguous diffeology and an unambiguous topology.
Indeed, we can first take the subspace differential structure on $Y$
and then the diffeology on $Y$ that it induces,
or we can first take the diffeology on $X$ that $\calF$ induces
and then the subset diffeology on~$Y$;
these two procedures yield the same diffeology on $Y$.
Also, 
the initial topology corresponding to the subspace differential structure on $Y$
coincides with the subspace topology on $Y$
induced by the initial topology on $X$. 
\eor
\end{proposition}

Diffeologies are not as well adapted to subsets:

\begin{remark}
\labell{r:subsets diffeological}
Given a diffeological space $(X,\calD)$ and a subset $Y \subset X$,
the subset $Y$ might not acquire an unambiguous differential structure
nor an unambiguous topology.
The two procedures --- first passing 
to the induced differential structure on $X$
and then to the subspace differential structure on $Y$, 
or first passing to the subset diffeology on $Y$
and then to the induced differential structure on $Y$ --- might 
yield two different differential structures on~$Y$.
Also, the D-topology corresponding to the subset diffeology on $Y$
might differ from the subset topology induced by the D-topology on $X$.
Both of these ambiguities occur with the subset $\Q$ of $\RR$
of Example \ref{e:rational numbers},
as well as with the ``pinched topologist's sine curve''
of Example~\ref{x:top sine curve}. 
\eor
\end{remark}

\begin{definition}[Quotients]\labell{d:quotients}
Let $X$ be a set, let $\sim$ be an equivalence relation on $X$,
and let $\pi \colon X \to X/\!\!\sim$
be the quotient map.
Given a differential structure $\calF$ on $X$,
the \emph{quotient differential structure} on $X/\!\!\sim$
consists of those functions $f \colon X/{\sim} \to \RR$
whose pullback $f \circ \pi \colon X \to \RR$ is in $\calF$.
Given a diffeology $\calD$ on $X$,
the \emph{quotient diffeology} on $X/\!\!\sim$ 
consists of those parametrisations $p \colon U \to X/\!\!\sim$ 
that locally lift to $X$ in the following sense: 
for every $u \in U$ there exist an open neighbourhood $V$ of $u$ in $U$
and a plot $q \colon V \to X$ such that $p|_V = \pi \circ q$.
\eod
\end{definition}

Diffeologies are well adapted to quotients:

\begin{proposition}[Diffeological quotients]\labell{p:quotients}
Given a diffeological space $(X,\calD)$ 
and an equivalence relation $\sim$ on $X$,
we obtain on the quotient $X/{\sim}$ 
an unambiguous differential structure and an unambiguous topology.
Indeed, we can first take the quotient diffeology on $X/{\sim}$
and then the differential structure that it induces,
or we can first take the differential structure on $X$ that $\calD$ induces
and then take the quotient differential structure on $X/{\sim}$.
These two procedures yield the same differential structure on $X/{\sim}$.
Also, the D-topology corresponding to the diffeology on $X/{\sim}$
coincides with the quotient topology on $X/{\sim}$ 
induced by the D-topology on $X$.
\eor
\end{proposition}

Differential structures are not as well adapted to quotients:

\begin{remark} \labell{r:quotients differential}
Given a differential space $(X,\calF)$ 
and an equivalence relation $\sim$ on $X$,
the quotient $X/{\sim}$ might not acquire 
an unambiguous diffeology nor an unambiguous topology.
The two procedures --- first passing 
to the induced diffeology on $X$
and then to the quotient diffeology on $X/{\sim}$,
or first passing to the quotient differential structure on $X/{\sim}$
and then to the induced diffeology on $X/{\sim}$ --- might 
yield two different diffeologies on the quotient $X/{\sim}$.
For example, this occurs with the irrational torus $\R /(\Z + \alpha \Z)$ 
as in Example~\ref{x:irrational flow},
with the quotient $\R/\Z_2$ as in Example~\ref{x:orbifolds},
and with the quotient $\R/(0,1)$ of the real line $\R$ 
by the open interval $(0,1)$ as in Example~\ref{x:R modulo open}.
Also, the initial topology corresponding to 
the quotient differential structure on $X/{\sim}$
might differ from the quotient topology on $X/{\sim}$
induced by the initial topology on $X$;
for example, this occurs with the quotient $\R/(0,1)$.
\eor
\end{remark}

We conclude this section with a couple of examples and open questions.

We start with an important collection of sub-quotients:

\begin{example}[Reduced Spaces]\labell{x:reduced spaces}
For a symplectic manifold $(M,\omega)$ with an action of a compact Lie
group $G$ and momentum map $\mu\colon M\to\mathfrak{g}^*$, the reduced
space $\mu^{-1}(0)/G$ inherits from $M$ an unambiguous diffeology
and an unambiguous differential structure, which are compatible.  
The differential structure on the reduced space $\mu^{-1}(0)/G$ 
does not always induce the diffeology
on the reduced space $\mu^{-1}(0)/G$ (Example~\ref{x:orbifolds}). 
\eor
\end{example}

Example~\ref{x:reduced spaces} raises an interesting question:

\begin{question}
In the setup of Example~\ref{x:reduced spaces},
does the diffeology on the reduced space $\mu^{-1}(0)/G$ 
necessarily induce the differential structure on the reduced space
$\mu^{-1}(0)/G$?
\end{question}

The following example illustrates
that diffeology can carry rich information about quotients
and that differential structures can carry rich information about subsets.

\begin{example}[Irrational flow on the torus]
\labell{x:irrational flow}
Fix an irrational number $\alpha$.
Consider the linear flow with slope $\alpha$ on the torus $\RR^2/\ZZ^2$:
$$ [x,y] \mapsto [x+t,y+\alpha t].$$
Let $T_\alpha$ be the quotient of the torus by this linear flow,
equipped with the quotient diffeology
(which induces the quotient differential structure 
and the quotient topology; see Proposition~\ref{p:quotients}).
Let $L_\alpha$ be the orbit through $[0,0]$ 
of this linear flow, equipped with the subspace differential structure
(which induces the subset diffeology and the subset differential structure;
see Proposition~\ref{p:subsets}).

The differential structure on $T_\alpha$ is trivial: 
it consists of the constant functions.
In contrast, the diffeology of $T_\alpha$ is non-trivial. 
For example, 
$$ t \mapsto \begin{cases}
 [0,0] & t<0 \\ [0,r] & t \geq 0
\end{cases} $$
is not a plot of $T_\alpha$ if $r \not\in \Z + \alpha\Z$.
Thus, the diffeological space $T_\alpha$ is not reflexive.

The diffeology on $L_\alpha$ is standard:
the inclusion map $t \mapsto [t,\alpha t]$ is a diffeomorphism 
from the real line $\RR$ with its standard diffeology to $L_\alpha$.
In contrast, the differential structure on $L_\alpha$ is not standard:
for example, its topology is not locally connected.
It follows that the differential space $L_\alpha$ is not reflexive.
\eoe
\end{example}

We now elaborate on Example~\ref{x:irrational flow},
leading to an open question.

An automorphism of the torus (as a Lie group)
carries the linear flow with slope $\alpha$
to a linear flow with slope $\beta$ where $\beta$ is obtained from $\alpha$ 
by a fractional linear transformation with integer coefficients:
$$ \beta = \frac{a \alpha + b}{c \alpha + d} \quad , \quad 
   a,b,c,d \in \ZZ \quad , \quad ad-bc = \pm 1 .$$
When $\alpha$ and $\beta$ are related in this way,
we say that they are $\GL(2,\ZZ)$-congruent.
Thus, if $\alpha$ and $\beta$ are $\GL(2,\ZZ)$-congruent, 
then the quotients $T_\alpha$ and $T_\beta$ are diffeomorphic
as diffeological spaces (hence also as differential spaces),
and the subsets $L_\alpha$ and $L_\beta$ are diffeomorphic 
as differential spaces (hence also as diffeological spaces).
Donato and Iglesias \cite{donato-iglesias} proved a striking result:
if $T_\alpha$ and $T_\beta$ are diffeomorphic as diffeological spaces, 
then $\alpha$ and $\beta$ are $\GL(2,\ZZ)$-congruent.
See Iglesias's book \cite[Exercise~4 with solution at the back of the book]
{iglesias}.

\begin{question}
Assuming that 
$L_\alpha$  and $L_\beta$ are diffeomorphic as differential spaces,
can we conclude that $\alpha$ and $\beta$ are $\GL(2,\ZZ)$-congruent?
\end{question}


\section{Orbifolds, Quotients by Compact Group Actions, 
and Manifolds with Corners}
\labell{sec:manifolds}

In this section, we study the quotient diffeological and quotient differential structures on the orbit space of a linear action of a compact Lie group on a cartesian space, and apply this to proper Lie group actions, to orbifolds, and to manifolds-with-corners.

\begin{example}[Orthogonal quotient] \labell{x:orthogonal action}
Let $G$ be a compact Lie group acting linearly on~$\RR^n$.
By a theorem of Hilbert \cite[p.~618]{weyl:hilbert},
the ring of $G$-invariant polynomials on $\RR^n$ is finitely-generated.
A choice of $m$ generators for this ring induces a $G$-invariant proper map
$i \colon \RR^n \to\RR^m$, which we call a \emph{Hilbert map}.
By a theorem of Gerald Schwarz \cite{schwarz}, 
every $G$-invariant smooth function on $\RR^n$
can be expressed as the pullback by $i$ of a smooth function on $\RR^m$.
This implies that the Hilbert map 
descends to a diffeomorphism from $\RR^n/G$, 
with the quotient differential structure induced from $\RR^n$,
to the image of the Hilbert map, 
with the subspace differential structure induced from $\RR^m$.

The quotient differential structure on $\RR^n/G$
is induced by the quotient diffeology on $\R^n/G$ 
by Proposition~\ref{p:quotients},
so it is reflexive by Proposition \ref{p:reflexive}.
Consequently, the subspace differential structure on the image of the
Hilbert map is reflexive.

In contrast, the quotient diffeology on $\RR^n/G$ 
might not be reflexive.  
For example, the map
$\RR^n /\O(n) \to [0,\infty)$ given by $[x] \mapsto \| x \|^2$
is an isomorphism of differential spaces (by Schwarz's theorem),
but the quotient diffeologies on $\RR^n / \O(n)$ are non-isomorphic
for different values of $n$ (see \cite{iglesias},
Exercise 50, with solution at the back of the book).
In particular, 
this Hilbert map does not induce a diffeomorphism of diffeological spaces 
from $\R^n/G$ to its image in $\R^m$.  
\eoe
\end{example}

\begin{example}[Proper Lie group action]\labell{x:cpct gp action} \ 
Combining Proposition~\ref{p:quotients}, Example~\ref{x:orthogonal action},
and the slice theorem \cite{koszul,palais},
the differential structure on the quotient 
of a manifold by a compact (or proper) Lie group action
is reflexive and is subcartesian; i.e., 
locally diffeomorphic to subsets of cartesian spaces. 
%
%
\eoe
\end{example}

\begin{example}[$\ZZ_2$- and $(\ZZ_2)^n$-actions]\labell{x:schwarz}
We note two special cases of Schwarz's theorem \cite{schwarz},
which in the case $n=1$ were proved by Whitney \cite{whitney:even}.
\begin{enumerate}
\item
Let the two-element group $\ZZ_2$ act on $\RR^n$ by
$(x_1,\ldots,x_n) \mapsto \pm (x_1,\ldots,x_n)$.
Then every invariant smooth function has the form
$g((x_ix_j)_{1 \leq i \leq j \leq n})$
where $g \colon \RR^{n(n+1)/2} \to \RR$ is smooth.
Here the Hilbert map $\RR^n \to \R^{n(n+1)/2}$
is given by
$(x_1,\ldots,x_n) \mapsto ( (x_ix_j)_{1 \leq i \leq j \leq n} )$.
When $n=2$, after a linear change of coordinates,
the image of the Hilbert map becomes the subset
$\{ z^2 = x^2 + y^2, \ z \geq 0\}$ of $\RR^3$.

\item
Let $(\ZZ_2)^n$ act on $\RR^n$ by
$(x_1,\ldots,x_n) \mapsto (\pm x_1, \ldots, \pm x_n)$.
Then every invariant smooth function has the form $g(x_1^2,\ldots,x_n^2)$
where $g \colon \RR^n \to \RR$ is smooth.
Here the Hilbert map $\RR^n \to \RR^n$ is given by 
$(x_1,\ldots,x_n) \mapsto (x_1^2,\ldots,x_n^2)$.
Its image is the positive orthant, $\RR_{\geq 0}^n$.
\eoe
\end{enumerate}
\end{example}

\begin{example}[Orbifolds] \labell{x:orbifolds}
(Effective) orbifolds can be defined
as diffeological spaces that are locally diffeomorphic
to quotients of the form $\RR^n/\Gamma$, where $\Gamma$
is a finite subgroup of $\O(n)$ (see \cite{IZKZ}).
As a differential space, an orbifold is reflexive.
However, the diffeology on an orbifold is generally not reflexive,
as illustrated in the following two examples.

Let the two-element group $\ZZ_2$ act on $\RR$ by $x \mapsto \pm x$,
and let $\pi \colon \R \to \R/\Z_2$ be the quotient map.
The quotient diffeology $\calD_{\R/\Z_2}$ induces the quotient
differential structure $\CIN(\RR/\ZZ_2)$
(see Proposition~\ref{p:quotients}),
but it is not induced by this differential structure:
the map $p(u,v) := [\pm\sqrt{u^2+v^2}]$ from $\RR^2$ to $\RR/\ZZ_2$ 
does not have a smooth lift near the origin,
but it is in the diffeology that is induced by $\CIN(\RR/\ZZ_2)$.
Indeed, by Schwarz's theorem (see Example~\ref{x:schwarz}),
if $f \in \CIN(\RR/\ZZ_2)$,
then $\pi^* f(x) = g(x^2)$ for some smooth function $g \colon \RR \to \RR$,
and so $f \circ p$ is equal to $(u,v) \mapsto g(u^2+v^2)$, which is smooth.
This shows that $\calD_{\RR/\ZZ_2}$ is not reflexive.

The following example is due to Moshe Zadka.
Let the two-element group $\ZZ_2$ act on $\RR^2$ by $(x,y) \mapsto \pm (x,y)$,
and let $\pi \colon \RR^2 \to \RR^2 / \ZZ_2 $
be the quotient map.
The quotient diffeology $\calD_{\RR^2/\ZZ_2}$ 
induces the quotient differential structure $\CIN(\RR^2/\ZZ_2)$
(see Proposition~\ref{p:quotients}),
but it is not induced by this differential structure:
the map 
$$ p (r\cos\theta,r\sin\theta) := \begin{cases}
[e^{-1/r^2}\cos(\theta/2),e^{-1/r^2}\sin(\theta/2)] & r>0 \\
[0,0] & r=0
\end{cases} $$
from $\RR^2$ to $\RR^2/\ZZ_2$
does not have a smooth (nor even continuous) lift near the origin,
but it is in the diffeology that is induced by $\CIN(\RR^2/\ZZ_2)$.
Indeed, by Schwarz's theorem (see Example \ref{x:schwarz}), 
if $f \in \CIN(\RR^2/\ZZ_2)$,
then $\pi^*f(x,y) = g(x^2,xy,y^2)$ for some smooth function 
$g \colon \RR^3 \to \RR$,
and so $f\circ p$ is equal to
$$ (r\cos\theta,r\sin\theta) \mapsto
\begin{cases}
g( e^{-2/r^2} \frac{1+\cos\theta}{2} , \
   e^{-2/r^2} \frac{\sin\theta}{2} , \
   e^{-2/r^2}\frac{1-\cos\theta}{2}) & r>0 \\
g(0,0,0) & r=0 ,
\end{cases}$$
which is smooth. 
This shows that $\calD_{\RR^2/\ZZ_2}$ is not reflexive.
\eoe
\end{example}

\begin{example}[The positive orthant] \labell{x:orthant}
On the positive orthant $\RR_{\geq 0}^n$,
the subspace differential structure that is induced from $\RR^n$ is reflexive.
Indeed, by Example~\ref{x:schwarz},
the positive orthant is the image of a Hilbert embedding,
and by Example~\ref{x:orthogonal action}, 
this implies that the differential structure $\calF$ is reflexive.
\eoe
\end{example}

\begin{example}[Manifolds-with-corners] \labell{x:corners}
We recall the definition of a manifold-with- corners.
An $n$-dimensional \emph{chart-with-corners} on a topological space $M$
is a homeomorphism $\varphi \colon U \to \Omega$
from an open subset $U$ of $M$ to a relatively open subset $\Omega$
of the positive orthant $\RR_{\geq 0}^n$.
Charts-with-corners $\varphi_1,\varphi_2$ are \emph{compatible}
if $\varphi_2 \circ \varphi_1^{-1}$
and $\varphi_1 \circ \varphi_2^{-1}$,
which are homeomorphisms between relatively open subsets of $\RR_{\geq 0}^n$,
are smooth in the sense that they locally extend
to smooth functions from $\RR^n$ to $\RR^n$.
An \emph{atlas-with-corners} on $M$ 
is a set of pairwise compatible charts with corners
whose domains cover $M$.  A \emph{manifold-with-corners}
is a Hausdorff, second-countable topological space $M$
equipped with a maximal atlas with corners.

An equivalent definition of manifold-with-corners is as a
Hausdorff, second countable differential space
that is locally functionally diffeomorphic to open subsets
of $\RR_{\geq 0}^n$.
A map between manifolds-with-corners is smooth in the classical sense
if and only if it is functionally smooth.
The differential structure on a manifold-with-corners is reflexive;
this follows from Example~\ref{x:orthant} (the positive orthant)
and the existence of smooth bump functions.

Manifolds-with-corners can also be viewed as diffeological spaces.
The D-topology on the positive orthant $\RR_{\geq 0}^n$
coincides with the subspace topology induced from $\R^n$;
this follows from the fact that the plot
$(x_1,\ldots,x_n)\mapsto (x_1^2,\ldots,x_n^2)$
restricts to a homeomorphism from the positive orthant to itself
with respect to the subspace topology.
By Proposition~\ref{p:subsets} and Example~\ref{x:orthant},
the subset diffeology and the subspace differential structure
on the positive orthant $\RR_{\geq 0}^n$ induce each other.
It follows that a map between relatively open subsets of the positive orthant
is a diffeological diffeomorphism if and only if it is a functional 
diffeomorphism,
which is equivalent to being a diffeomorphism in the classical sense.
It further follows that a manifold-with-corners can be equivalently defined 
as a diffeological space
that is locally diffeomorphic to open subsets of $\RR_{\geq 0}^n$.
\eoe
\end{example}

\smallskip\noindent\textbf{Notes.}

\begin{enumerate}

\item
The argument in Example~\ref{x:orthant} is a generalisation of the same
statement for half-spaces $\RR^{n-1}\times[0,\infty)$
that was given by Iglesias-Zemmour in \cite[ch.~4]{iglesias} to show that
the classical notion of a manifold-with-boundary is the same 
as the diffeological notion.
This generalisation appeared in \cite{LW}.  
More generally, 
the subspace differential structure on any locally closed convex set 
is reflexive;
see \cite{KW:convex}.

\item
Manifolds-with-corners were introduced in 1961 by Jean Cerf 
and by Adrien Douady \cite{cerf,douady}
and are now included in standard textbooks
such as John Lee's \cite[Chapter~16]{lee}.
Our definitions of a manifold-with-corners are equivalent to theirs.
These definitions are local.

Manifolds-with-faces were introduced in 1968 by Klaus J\"anich \cite{janich};
also see \cite[Chap.~4]{verona}.
The codimension-$k$ \emph{strata} of a manifold-with-corners
are the connected components of the set of those points
that, in a chart-with-corners, have exactly $k$ coordinates that vanish.
Manifolds-with-faces are manifolds-with-corners 
in which every codimension-$k$ stratum 
is in the closure of $k$ distinct codimension one strata.
This condition is global.
The disc-with-one-corner in the plane,
given in polar coordinates by $r \leq \sin 2\theta$
for $0 \leq \theta \leq \pi/2$,
is a manifold-with-corners but is not a manifold-with-faces.
Some authors
use the term ``t-manifold'' for a manifold-with-corners
and the term ``manifold-with-corners'' for a manifold-with-faces;
see \cite[Definition 1.8.5]{melrose};
also see \cite[Article 1.1.19]{pflaum}.

\end{enumerate}

\section{Intersecting submanifolds}
\labell{sec:wedge}

In this section we consider the diffeological and differential structures on some unions of lines in the plane, and more generally, on some unions of submanifolds in an ambient manifold.

\begin{example}[Two coordinate axes] \labell{x:wedgesum}
Consider the wedge sum of two copies of $\RR$ attached at their origins,
which we write as 
$ X = ( \RR_1 \amalg \RR_2) / (0_1 \sim 0_2) ;$
denote the quotient diffeology by $\calD_X$.
Let $E \subset \R^2$ be the union of the two coordinate axes 
in the cartesian plane,
equipped with its subspace differential structure.  
Let
$$ \varphi \colon X \to E $$
be the bijection 
whose pullback to $\RR_1$ is $x \mapsto (x,0)$
and whose pullback to $\RR_2$ is $y \mapsto (0,y)$.
Then
\begin{enumerate}
\item
The map $\varphi$ is a functional diffeomorphism
from the differential space $(X,\Phi\calD_X)$ to the differential subspace $(E,\CIN(E))$ of $\R^2$.
Moreover, the differential structure $\Phi\calD_X$
consists of those real-valued functions on $X$
whose pullbacks to $\R_1$ and to $\R_2$ are smooth.
\item
The differential space $(E,\CIN(E))$ is reflexive.
\item
The diffeological space $(X,\calD_X)$ is not reflexive.
\eoe
\end{enumerate}
\end{example} 

\begin{proof}
We prove Item (1) in \S\ref{subsec:lines in plane}.
By Item (1), $\CIN(E)$ is a differential structure
that is induced by some diffeology;
Proposition~\ref{p:reflexive} (``reflexive stability'')
then gives Item (2).
For Item (3) we need to show that $\Pi\Phi\calD_X \supsetneq \calD_X$.
Consider 
the parametrisation $p \colon \RR \to X$ whose composition with $\varphi$
is
$$ t \mapsto \begin{cases}
 ( e^{-1/t^2} , 0 ) & \text{ if } t < 0 \\
 ( 0 , 0 ) & \text{ if } t = 0 \\
 ( 0 , e^{-1/t^2} ) & \text{ if } t > 0 .
\end{cases} $$
Because this composition is a smooth map with image in $E$,
it is a plot of $E$; by Proposition~\ref{p:subsets} it is in $\Pi\CIN(E)$;
Item~(1) implies that $p$ is in $\Pi\Phi\calD_X$.
On the other hand, $p$ does not lift
to a smooth (nor even continuous) map to $\RR_1 \amalg \RR_2$
on any neighbourhood of $t=0$,
so $p$ is not in the quotient diffeology $\calD_X$ on $X$. This proves~(3).
\end{proof}

\begin{remark}\labell{r:wedgesum}
\noindent
\begin{enumerate}
	\item Example~\ref{x:wedgesum} generalises to any finite number of copies of $\RR$.  In particular, the subspace differential structure on the union of the three coordinate axes in $\RR^3$ is reflexive.

	\item Example~\ref{x:wedgesum} generalises to arbitrary pointed manifolds
$(N_1,*_1) , \ldots , (N_k,*_k)$, of possibly different dimensions,
with $X := (N_1 \amalg \ldots \amalg N_k) / *_i \sim *_j$
for all $i,j$, and with $E \subset N_1 \times \ldots \times N_k$.
See \cite[Example 2.67]{watts}.

	\item Example~\ref{x:wedgesum}
also generalises to transversally intersecting submanifolds;
see \cite[Example 2.70]{watts}.
Here, we take $N_1$ and $N_2$ to be (embedded) submanifolds 
of an ambient manifold $M$.
We let $i \colon N_1 \amalg N_2 \to M$ be the map whose restriction
to each $N_i$ is the inclusion map, we take $X := N_1 \amalg N_2 /{\sim}$
where $x \sim y$ if and only if $i(x) = i(y)$,
and we take $E := N_1 \cup N_2 \subset M$.
\eor
\end{enumerate}
\end{remark}

The generalisations mentioned in Remark~\ref{r:wedgesum} are special cases of the following more general example:

\begin{example}[Cleanly-intersecting submanifolds]\labell{x:wedgesum2}
\noindent
Let $M$ be a manifold,
and let $\{ i_\tau \colon N_\tau \hookrightarrow M \}$
be a family of submanifolds 
whose intersections are jointly clean in the following sense.
Each point of $M$ is in the domain of some coordinate chart
$ \varphi \colon U \to \Omega \subseteq \R^n$
such that, for each $\tau$, if the submanifold $N_\tau$ meets $U$,
then $\varphi(U \cap N_\tau)$ is the intersection of $\Omega$
with a coordinate subspace of $\R^n$.
Let $ X := \left( \coprod\limits_{\tau} N_\tau \right)/{\sim} $
where, for $x \in N_{\tau}$ and $y \in N_{\tau'}$,
we have $x \sim y$ if and only if $i_\tau(x)=i_{\tau'}(y)$;
denote the quotient diffeology by $\calD_X$.
Consider the subset $ E := \bigcup\limits_{\tau} i_\tau(N_\tau) $ of $M$,
let $\CIN(E)$ be the subspace differential structure,
and let $ \varphi \colon X \to E $
be the bijection whose pullback to $N_\tau$ is $i_\tau$.
Then, as in Example~\ref{x:wedgesum},
\begin{enumerate}
\item
The map $\varphi$ is a functional diffeomorphism
from the differential space $(X,\Phi\calD_X)$
to the differential subspace $(E,\CIN(E))$ of $M$.
Moreover, the differential structure $\Phi\calD_X$
consists of those real-valued functions on $X$
whose pullback to each $N_\tau$ is smooth.
\item
The differential space $(E,\CIN(E))$ is reflexive.
\item
The diffeological space $(X,\calD_X)$ is not reflexive,
unless the components of $E$ are submanifolds of~$M$.
\end{enumerate}
For details, see \S\ref{subsec:lines in plane}.
\end{example}

\begin{example} [Three lines in $\RR^2$] \labell{three lines}
Let $S$ be the subset of $\RR^2$ given by the union of the $x$-axis, 
the $y$-axis, and the line $y = x$,
with the subspace differential structure $\CIN(S)$
and the subset diffeology $\calD_S$ that are induced from $\RR^2$.
Let $E \subseteq \RR^3$ be the union of the three coordinate axes,
with the subspace differential structure $\CIN(E)$
and the subset diffeology $\calD_E$
that are induced from $\RR^3$.
Consider the bijection
$ \varphi \colon E \to S $
given by $(t,0,0) \mapsto (t,0)$, \ 
$(0,t,0) \mapsto (0,t)$, and $(0,0,t) \mapsto (t,t)$.  Then 
\begin{enumerate}[noitemsep,topsep=0pt]
\item \labell{partA}
The map $\varphi$ is a diffeological diffeomorphism
from $(E,\calD_E)$, where $\calD_E$ is the subset diffeology on $E$
that is induced from $\R^3$,
to $(S,\calD_S)$, where $\calD_S$ is the subset diffeology on $S$
that is induced from $\R^2$.

\item \labell{partB}
The differential space $(S,\CIN(S))$ is not reflexive.
\eoe
\end{enumerate}
\end{example}

\begin{proof}
Because $\varphi$ extends to a smooth map between the ambient spaces
$\RR^3 \to \RR^2$,
(for example, take $(x,y,z) \mapsto (x+z,y+z)$,)
the map $\varphi \colon E \to S$ is functionally and diffeologically smooth,
so $\varphi \circ \calD_E \subseteq \calD_S$.
We prove the opposite inclusion,
$ \calD_S \subseteq \varphi \circ \calD_E $,
in \S\ref{subsec:lines in plane}.

By Item 1 of Remark~\ref{r:wedgesum}, the differential structure $\CIN(E)$ on $E$ is reflexive.
We then have
\[ \begin{array}{lll}
(S,\Phi\Pi\CIN(S))  & = (S,\Phi \calD_S) 
            & \text{ since } \calD_S = \Pi\CIN(S)
                       \text{ (see Proposition~\ref{p:subsets})}, \\
  & \cong (E,\Phi\calD_E) & \text{ by Part~\eqref{partA}}, \\ 
  & = (E,\CIN(E)) & \text{ since } \calD_E  = \Pi\CIN(E)
                     \text{ and } \CIN(E) \text{ is reflexive. } 
\end{array} \]
The dimension of the Zariski tangent space at a point
in a differential space
is invariant under functional diffeomorphisms (see \cite{LSW}).
Since the dimension of the Zariski tangent space 
at the origin in $S$ is $2$, and that at the origin in $(E,\CIN(E))$, 
hence in $(S,\Phi\Pi\CIN(S))$, is $3$, 
the differential space $(S,\CIN(S))$ is not reflexive.
\end{proof}

\begin{example}[Many lines in $\RR^2$] \labell{x:lines in Rtwo}
For any integer $k \geq 3$,
Example~\ref{three lines}
generalises to the union of any $k$ distinct lines through the 
origin in $\R^2$.
In particular,
every two such unions are diffeomorphic as diffeological spaces.
In contrast, 
two such unions are diffeomorphic as differential spaces
if and only if they differ by a linear transformation of $\R^2$.
(Given a diffeomorphism between them, take its differential at the origin.)
Thus, for each $k \geq 4$, such unions produce a continuum of non-isomorphic
differential spaces.
\end{example}

{
The following example was communicated to us by Katrin Wehrheim
as a topological space that can arise in the context of polyfolds.
We're interested in its diffeology:

\begin{example}[Axis and half-plane] \labell{x:Katrin}
Take the space $X$ that is obtained by gluing the $x$-axis
$\{ (x,0) \ | \ x \in \R \}$ with the open right half-plane 
$\{ (x,y) \ | \ x > 0 , \ y \in \R \}$ along their intersection in $\R^2$.
By Proposition~\ref{p:quotients}, its quotient diffeology $\calD_X$
induces its quotient differential structure $\calF_X$
and its quotient topology.
These are not induced by the natural inclusion map $X \hookrightarrow \R^2$:
the function 
$$f(x,y) := \begin{cases} 0 & y=0 \\ \frac{e^{-1/|y|}}{x} & y\neq 0, \ x>0
\end{cases}$$
is in $\calF_X$ but it is not continuous with respect to the subset topology
induced from $\R^2$.
\end{example}

\begin{question} \labell{q:Katrin}
In Example~\ref{x:Katrin}, is the diffeology $\calD_X$ reflexive?
\end{question}

\smallskip\noindent\textbf{Notes.}
From the point of view of Fr\"olicher spaces, wedge products 
are also analyzed in Batubenge and Ntumba's paper \cite[pages 76--78]{NB}.
Other aspects of the above examples also appeared in 
Watts' Ph.D.\ thesis~\cite[Examples~2.67 and~2.70]{watts}
and in Christensen and Wu's paper 
\cite[Examples~3.17, 3.19, and 3.20]{CW:tangent}.

{

\section{Topological Considerations}
\labell{sec:topology}

Recall that the D-topology of a diffeological space
is the strongest topology making all of its plots continuous, 
and the initial topology of a differential space
is the weakest topology 
making all of the functions of its differential structure continuous.  
The purpose of this section is to point out some properties
of these topologies that are necessary for reflexivity.

We begin with a simple observation:

%
\begin{lemma}[Compatible topologies]
\labell{l:top contain}
Let $X$ be a set, and let $\calD$ and $\calF$ be a diffeology and a
differential structure on $X$, respectively.  Suppose that $\calD$ and
$\calF$ are compatible.  Then the initial topology induced by $\calF$
is contained in the D-topology induced by $\calD$.
\end{lemma}
\begin{proof}
Let $f \in \calF$ and $p \in \calD$. Because $f \circ p$ is 
(smooth, hence) continuous, for any open interval $I$
the preimage $p^{-1}(f^{-1}(I))$ is open in the domain of $p$.
Since $p \in \calD$ is arbitrary, $f^{-1}(I)$ is D-open.
Since $f$ is arbitrary, the D-topology contains the initial topology.
\end{proof}

Propositions~\ref{p:d-top locally-path-conn}
and~\ref{p:init-top smoothly reg} below are known \cite{hector,laubinger};
for completeness, we include their proofs.

%
\begin{proposition}[Diffeological spaces are locally path-connected]
\labell{p:d-top locally-path-conn}
The D-topology of a diffeological space is locally path-connected.  
(Consequently, the connected components coincide with the path-connected 
components, and these components are both open and closed.)
\end{proposition}

\begin{proof}
Let $X$ be a diffeological space, let $x$ be a point of $X$,
and let $V$ be a D-open neighbourhood of $x$.
We need to show that  $V$  contains a path-connected neighbourhood of $x$.
Let  $X'$ be the smooth path component of  $x$  in $V$;
we will show that  $X'$  is  D-open.
Let  $p \colon U \to V$  be a plot.  For each connected component  $U'$  of  $U$,
the image  $p(U')$  is smoothly path-connected, so it is either contained in  $X'$
or disjoint from  $X'$.   So the preimage  $p^{-1}(X')$, being a union
of connected components of  $U$,  is open in  $U$.
Varying the plot  $p$,  this shows that the preimage  $X'$  is D-open.
Since $X'$ is smoothly path-connected, it is path-connected.
\end{proof}

Recall that a topological space $X$ is \textbf{completely regular}
if for every closed set $C$ and point $x \in X \ssminus C$
there exists a continuous function $f\colon X\to[0,1]$ 
that vanishes on $C$ and is equal to $1$ on $x$.
Spaces that are $T_0$ (points are distinguishable by open sets)
and completely regular are called $T_{3\frac{1}{2}}$.
Such spaces are automatically $T_3$ (regular and $T_0$),
hence $T_2$ (Hausdorff), hence $T_1$ (points are closed).
%

\begin{proposition}[Differential spaces are completely regular]
\labell{p:init-top smoothly reg}
The initial topology of a differential space is completely regular.
Consequently, any $T_0$ differential space is Hausdorff.
\end{proposition}

\begin{proof}
Let $(X,\calF)$ be a differential space, equipped with the initial topology.
Let $C$ be a closed subset of $X$, and let $x$ be a point in $X \ssminus C$.
By the definition of the initial topology,
there exist functions $h_1,\ldots,h_k \in \calF$
and open intervals $I_1,\ldots,I_k$
such that $x \in \cap_{i=1}^k h_i^{-1}(I_i) \subset X \ssminus C$.
Take $f := b \circ (h_1,\ldots,h_k)$
where $b \colon \R^k \to [0,1]$ is a smooth function
whose support is contained in $I_1 \times \ldots \times I_k$
and such that $f(h_1(x),\ldots,h_k(x)) = 1$.
\end{proof}

Lemma~\ref{l:top contain} and Propositions~\ref{p:d-top locally-path-conn} 
and~\ref{p:init-top smoothly reg} imply the following 
topological necessary conditions for the compatibility
of a diffeology and a differential structure.
%

\begin{corollary}\labell{c:top contain}
Let $X$ be a set, and let $\calD$ and $\calF$ be a diffeology 
and differential structure on $X$, respectively.  
Suppose that $\calD$ and $\calF$ are compatible.  
\begin{enumerate}
\item If the initial topology induced by $\calF$ is $T_0$, 
then the D-topology induced by $\calD$ is Hausdorff.
\item If the D-topology induced by $\calD$ is connected, 
then the initial topology induced by $\calF$ is path-connected.
\end{enumerate}
\end{corollary}

Using the ideas developed above, we show some necessary conditions
for reflexivity of diffeological spaces and of differential spaces.
%

\begin{proposition}[$T_0$ reflexive diffeological spaces]\labell{p:top obstr to refl diffeol}
Every $T_0$ reflexive diffeological space is Hausdorff.
\end{proposition}

\begin{proof}
Let $(X,\calD)$ be a diffeological space whose D-topology is $T_0$.
Let $x,y\in X$ be distinct points such that for any open neighbourhoods $U$
of $x$ and $V$ of $y$ we have $U\cap V\neq\emptyset$.
Since the D-topology is $T_0$, without loss of generality 
there exists an open neighbourhood $W$ of $y$ so that $x\notin W$.  
Define $p\colon\RR\to\{x,y\}$ by
$$p(t):=\begin{cases}
x & \text{~if $t<0$,}\\
y & \text{~if $t\geq 0$.}\\
\end{cases}$$
Then $p^{-1}(W)=[0,\infty)$.  
Since $p$ is not continuous, $p\notin\calD$.

Let $f\in\Phi\calD$.  Then $f(x)=f(y)$.  (Otherwise, 
setting $a=f(x)$, $b=f(y)$, and $0<\eps<\frac{|b-a|}{2}$, 
we have $x\in U:=f^{-1}((a-\eps,a+\eps))$ 
and $y\in V:=f^{-1}((b-\eps,b+\eps))$, 
and the intersection $U\cap V$ is empty.
By Lemma~\ref{l:top contain}, $U$ and $V$ are D-open.
This contradicts the choice of $x$ and $y$.)
So the composition $f\circ p$ is (constant, hence) smooth.
This shows that $p\in\Pi\Phi\calD$.
So $\calD$ is not reflexive.
\end{proof}
%

\begin{proposition}[Locally smoothly path-connected differential spaces]\labell{p:top obstr to refl diff str}
On every reflexive differential space,
the initial topology is locally smoothly path-connected.
\end{proposition}

\begin{proof}
Let $(X,\calF)$ be a reflexive differential space,
equipped with the initial topology, let $x \in X$ be any point,
and let $U$ be an open neighbourhood of $x$ in $X$.
Let $C$ be the smooth path component of $x$ in $U$.
It is enough to show that $C$ is a neighbourhood of $x$.

Let $b \in \calF$ be a smooth function such that $b(x) = 1$
and whose support $\supp(b)$ is contained in $U$
(cf.\ the proof of Proposition~\ref{p:init-top smoothly reg}).
Define $g \colon X \to \R$ to be equal to $b$ on $C$ and zero outside $C$.
Since $g^{-1}((0,\infty))$ contains 
$x$ and is contained in $C$,
it is enough to show that $g \in \calF$.

Let $p \in \Pi \calF$.  
Since the D-topology induced by $\Pi\calF$ contains the initial topology induced by $\calF$, the connected components of $p^{-1}(U)$, as well as $p^{-1}(X\smallsetminus\supp(b))$, are open in the domain of $p$.  Let $q$ be the restriction of $p$ to one of the connected components of $p^{-1}(U)$.  
If the image of $q$ does not meet $C$, then $g \circ q$ is identically $0$.  
Suppose the image of $q$ does meet $C$.
Since $C$ is a smooth path component of $U$, the image of $q$ is contained in $C$,
and so $g \circ q$ is equal to $b \circ q$, which is smooth.  Since $g\circ p$ is identically zero on the $p^{-1}(X\smallsetminus\supp(b))$, it is smooth.
It follows that $g \in \Phi\Pi \calF$.
Since $\calF$ is reflexive, $g \in \calF$.
\end{proof}

{

Here are a couple of applications
of the necessary conditions for reflexivity
that we gave in Propositions~\ref{p:top obstr to refl diffeol}
and~\ref{p:top obstr to refl diff str}.
%

\begin{example}[Line with double origin]
\labell{e:double origin}
Glue two copies of the real line along the complement of the origin;
write the quotient space as
$X := (\R_1 \amalg \R_2)/{\sim}$. 
Consider its quotient diffeology.
Its D-topology, which coincides with the quotient topology
(see Proposition~\ref{p:quotients}), is $T_0$ but not Hausdorff.
So the diffeological space $X$ is not reflexive; 
see Proposition~\ref{p:top obstr to refl diffeol}.
\end{example}
%

\begin{example}[Pinched Topologist's Sine Curve]\labell{x:top sine curve}
Let $Y \subset \R^2$ be the image of the curve
$\gamma \colon [0,1] \to \R^2$ that is given by 
$$ \gamma(t) = \begin{cases}
 ( \, 0, \, 0 \,) & \text{ if } t=0 \\
 ( \, t, \, t \sin(1/t) \, ) & \text{ if } 0 < t \leq 1 ,
\end{cases}$$
equipped with the subspace topology and differential structure $\calF_Y$ 
induced from $\RR^2$ (see Proposition~\ref{p:subsets}).
The point $(0,0)$ does not have any neighbourhood 
that is smoothly path connected;
this follows from the fact that 
the curve $\gamma$ is not rectifiable.
It follows from Proposition~\ref{p:top obstr to refl diff str} that
$\calF_Y$ is not reflexive.  
\eoe
\end{example}

}

Examples~\ref{x:R modulo closed} and \ref{x:R modulo open} below
originally appeared in Jordan Watts' thesis 
\cite[Examples 2.74 and 2.76]{watts}.  
In both of these examples, a quotient space 
is obtained from $\RR$ by collapsing an interval to a point.
 

\begin{example}[$\RR$ modulo a closed interval]
\labell{x:R modulo closed}
Let $X:=\RR/[0,1]$, equipped with the quotient diffeology $\calD_X$,
the quotient differential structure $\calF_X$, and the quotient topology.
The initial topology induced by $\calF_X$ 
coincides with the quotient topology.
This follows from the fact that the function $f \colon X \to \RR$ 
whose pullback to $\RR$ is
$$ x \mapsto \begin{cases}
 - e^{-\frac{1}{|x|}} & \text{ if } x<0 \\
 0                   & \text{ if } x \in [0,1] \\
 e^{-\frac{1}{|x-1|} } & \text{ if } x > 1
\end{cases} $$
is in $\calF_X$ and is a homeomorphism with respect to the quotient
topology on $X$. 
The map $\varphi \colon X \to \RR$ whose pullback to $\RR$ is 
\begin{equation*}
 x \mapsto \begin{cases}
 x & \text{ if } x < 0 \\
 0 & \text{ if } x \in [0,1] \\
 x-1 & \text{ if } x > 1
\end{cases}
\end{equation*}
is a diffeomorphism between the differential spaces $(X,\calF_X)$ 
and $(\RR,\calF)$,
where $\calF$ is the set of those smooth functions
on $\RR$ whose derivatives of all positive orders vanish at $0$.
The quotient diffeology $\calD_X$ is not reflexive:
$\varphi^{-1} \colon \RR \to X$ is in 
($\Pi\calF_X$, which by Proposition~\ref{p:quotients} is) $\Pi\Phi\calD_X$,
but it does not have any smooth (or even continuous) lift to $\R$
in any neighbourhood of the origin.
\eoe
\end{example}

\begin{example}[$\RR$ modulo an open interval] \labell{x:R modulo open}
Let $Y:=\RR/(0,1)$, equipped with the quotient diffeology $\calD_Y$.
Let $\pi_Y \colon \R \to Y$ be the quotient map.
The one-point set $\pi_Y((0,1))$ is open with respect
to (the quotient topology, hence) the D-topology.
In fact, this topology is $T_0$.
But this topology is not Hausdorff:
the points $\pi_Y(0)$ and $\pi_Y(1)$ do not have disjoint neighbourhoods. 
By Proposition~\ref{p:top obstr to refl diffeol}, 
$\calD_Y$ is not reflexive.

Recall that the quotient diffeology $\calD_Y$
induces the quotient differential structure $\calF_Y$.
In contrast with the previous example, the initial topology
induced by $\calF_Y$ is strictly smaller than the quotient topology:
in the initial topology, the points $\pi_Y(0)$, $\pi_Y((0,1))$, and $\pi_Y(1)$ 
are topologically indistinguishable.
The corresponding Kolmogorov quotient (see Remark~\ref{r:cloning} below)
can be identified with $\RR/[0,1]$.
\eoe
\end{example}

Examples~\ref{x:R modulo closed} and \ref{x:R modulo open}
motivate the following general remark on indistinguishable points.

{
%
\begin{remark}[Indistinguishable points]\labell{r:cloning}
%
Given a differential space $X$, one can create another differential space $Y$
by creating ``clones'' of points of $X$, which are not distinguishable 
by smooth functions.  In fact, up to isomorphism, such ``cloning'' 
is the only way of obtaining a differential space whose initial topology 
is not Hausdorff.
More precisely, let $(X,\calF)$ be a differential space.  
Let $X_K$ be the quotient of $X$ by the equivalence relation where $x \sim x'$
iff $f(x)=f(x')$ for all $f \in \calF$.
Equip $X_K$ with the quotient differential structure
(see Definition~\ref{d:quotients}).
Up to isomorphism, $X$ is obtained from $X_K$ by ``cloning''.
The initial topology of $X_K$ induced by $\calF_K$ 
coincides with its quotient topology 
(contrast with Remark~\ref{r:quotients differential}).
It is $T_0$ (hence, by Proposition~\ref{p:init-top smoothly reg},
it is $T_{3\frac{1}{2}}$).
The $T_0$ differential space $X_K$, with the map $\pi \colon X \to X_K$, 
satisfies the following universal property:
\begin{equation} \labell{e:universality}
\begin{minipage}{5in}
If $Y$ is a $T_0$ differential space and $\varphi \colon X \to Y$ 
is a functionally smooth map, 
then there exists a unique functionally smooth map
$\varphi_K \colon X_K \to Y$ such that $\varphi = \varphi_K \circ \pi $.
\end{minipage}
\end{equation}
Topologically, $X_K$ coincides with the \emph{Kolmogorov quotient} of $X$,
which is the quotient by the equivalence relation where $x \sim x'$
if and only if each open neighbourhood of $x$ contains $x'$ and vice versa.
It also coincides with the \emph{Hausdorffification} of $X$
(whose construction for more general topological spaces
may require iterated quotients, and possibly transfinite recursion; 
see \cite{vanM}).
These satisfy universal properties similar to~\eqref{e:universality}
but with respect to continuous maps to $T_0$ spaces and to Hausdorff spaces,
respectively. 
\eor
\end{remark}
}

\begin{example}[$\RR$ Modulo $(x\sim 2x)$]\labell{x:R mod 2}
Take $X = \R/{\sim}$ where $x \sim y$ if and only if $y = 2^m x$
for some integer~$m$.  
Its quotient diffeology is non-trivial, 
but its differential structure consists of the constant functions,
so its quotient diffeology is not reflexive.
Another way to see this is to note that the D-topology
coincides with the quotient topology (Proposition~\ref{p:quotients}),
which is $T_0$, but not Hausdorff, 
and to apply Proposition~\ref{p:top obstr to refl diffeol}.
\eoe
\end{example}

We have seen several non-reflexive quotients of reflexive diffeological spaces:
the irrational torus $\R/(\Z+\alpha\Z)$ (Example~\ref{x:irrational flow}),
the orbifold $\R^2/\Z_2$ (Example~\ref{x:orbifolds}),
the quotients $\R/[0,1]$ and $\R/(0,1)$
(Examples~\ref{x:R modulo closed} and~\ref{x:R modulo open}),
and the quotient $\R/(x \sim 2x)$ (Example~\ref{x:R mod 2}).
In the first three of these, 
the initial topology coincides with the quotient topology;
in the last two, the initial topology is different from the quotient topology. 
The second and third of these are Hausdorff; the others are non-Hausdorff.

These examples raise the question of whether 
the initial topology coincides with the quotient topology
on a quotient differential space that is Hausdorff.
The answer is no:
the following example, inspired by the Moore-Niemytzki plane
(see \cite[Example 82]{counterexamples}, where it is called the Niemytzki tangent disk topology), exhibits a quotient differential
space whose initial topology is Hausdorff but is strictly smaller
than its quotient topology.

\begin{example}[A Moore-Niemytzki-like topology]\labell{x:niemytzki}
Let $H$ be the open upper half plane in $\RR^2$,
and let $\ol{H}$ be its closure in $\R^2$.
Equip $H$, $\overline{H}$, and the sets 
$$ C_x := H\cup\{(x,0)\}, $$
for $x\in\RR$, with the subspace differential structures 
and subset diffeologies that are induced from $\RR^2$.
Equip
$$ X :=\coprod_{x\in\RR} C_x \, $$
with the coproduct differential structure, denoted $\calF_X$,
and with the coproduct diffeology, denoted $\calD_X$.
Let $Y$ be the gluing of the components of $X$ along $H$,
equipped with the quotient differential structure $\calF_Y$
and the quotient diffeology $\calD_Y$.
Denote the inclusion maps of the $C_x$ into $X$ 
and the quotient map from $X$ to $Y$ as
$$ i_x \colon C_x \to X \quad \text{ and } \quad \pi \colon X \to Y.$$

There exists a unique bijection
$$\varphi\colon Y\to\overline{H} $$
such that, 
for each $x \in \R$, the composition 
$\varphi \circ \pi \circ i_x \colon C_x \to \ol{H}$
is the inclusion map of $C_x$ into $\ol{H}$. 
Because this inclusion map is smooth,
and because the components $C_x$ form an open covering of $X$,
it follows that $\varphi \circ \pi \colon X \to \ol{H}$ is smooth.
By the definition of the quotient differential structure, 
it follows that $\varphi \colon Y \to \ol{H}$ is smooth.

We claim:
\begin{equation} \labell{CxH}
\begin{minipage}{4in}
For any function $g \colon \ol{H} \to \R$, \\
if $g|_{C_x} \colon C_x \to \R$ is functionally smooth for all $x$, \\
then $g$ is functionally smooth.
\end{minipage}
\end{equation}
Indeed,
let $g \colon \ol{H} \to \R$, and suppose that $g|_{C_x} \colon C_x \to \R$
is smooth for all $x$. Then $g|_H \colon H \to \R$ is smooth.
Let $x \in \R$. Because $g|_{C_x}$ is smooth,
there is an open neighbourhood $U_x$ of $(x,0)$ in $\R^2$
and a smooth function $h_x \in \Cinf(U_x)$ 
that coincides with $g$ on the subset $C_x \cap U_x$.
Let $(x',0) \in U_x$. 
For all sufficiently large $n$, we have $(x',\frac{1}{n}) \in U_x$,
and so $h_x(x',\frac{1}{n}) = g(x',\frac{1}{n})$.
But $h_x(x',\frac{1}{n}) \xrightarrow[n \to \infty]{} h_x(x',0)$
because $h_x$ is smooth on $U_x$,
and $g(x',\frac{1}{n}) \xrightarrow[n \to \infty]{} g(x',0)$
because $g$ is smooth on $C_{x'}$,
so $h_x(x',0) = g(x',0)$.
Because $x'$ was arbitrary,
$h_x$ coincides with $g$ on all of $\ol{H} \cap U_x$.
Because $x$ was arbitrary, it follows that $g$ is smooth.

Because the inverse
$ \varphi^{-1} \colon \ol{H} \to Y $
satisfies 
$$\varphi^{-1}|_{C_x} = \pi \circ i_x \colon C_x \to Y \, $$
for all $x$,
it follows from \eqref{CxH} that this inverse is functionally smooth.
Thus, $\varphi$ is a functional diffeomorphism.

By Christensen-Sinnamon-Wu \cite[Lemma 3.17]{CSW},
since each $C_x$ is a convex subset of $\R^2$, 
its D-topology is equal to its subspace topology induced from $\R^2$.
Because the 
D-topology induced by the quotient diffeology $\calD_Y$ on $Y$
is equal to the quotient topology $\tau_Y$ on $Y$ induced from $X$
(Proposition~\ref{p:quotients}),
a subset $A$ of $Y$ is closed if and only if 
$i_x^{-1}(\pi^{-1}(A))$ is closed in $C_x$ for all $x$.
Consequently, every subset of $\varphi^{-1}(\partial H)$ is closed.
Thus, while $\varphi$ is a continuous bijection from $(Y,\tau_Y)$ to
$\overline{H}$, it is not a homeomorphism.

It follows that the initial topology on $Y$ 
is strictly smaller than the quotient topology $\tau_Y$.
Note, though, that both of these topologies are Haudorff.

Finally, note that the diffeology $\calD_Y$ is not reflexive.
Indeed, because $\varphi \colon Y \to \ol{H}$
is a functional diffeomorphism, 
the map $x \mapsto \varphi^{-1}(x,0)$ from $\R$ to $Y$ is functionally smooth.
But this map is not diffeologically smooth.
\eoe
\end{example}

\smallskip\noindent\textbf{Notes.}

\begin{itemize}
\item
In a diffeological space, the path components (with respect to the D-topology)
coincide with the (diffeologically) smooth path components
\cite[Article 5.7]{iglesias}.
\item
By the proof of Proposition~\ref{p:init-top smoothly reg}, 
every differential space is \emph{smoothly regular},
in the following sense~\cite{KM}:
for every closed set $C$ and point $x \in X \ssminus C$
there exists a smooth function $f \colon X \to [0,1]$ 
that vanishes on $C$ and is equal to $1$ on $x$.

\item
By \cite[Theorem 3.10]{virgin}:
given a diffeology $\calD$ and a differential structure $\calF$
that are compatible,
the D-topology coincides with the initial topology
if and only if the D-topology is smoothly regular.

\item
\jnote{added:}
In Example~\ref{x:niemytzki}, we showed that the D-topology of $Y$ is Hausdorff.  However, it is not completely regular.  Indeed, since any subset of $\varphi^{-1}(\del H)$ is closed, there is no continuous function that separates $\pi\circ i_x((x,0))$ from its complement in $\varphi^{-1}(\del H)$.  It follows that the D-topology of $Y$ is not smoothly regular.

\item
In Corollary~\ref{c:top contain} we can obtain
a stronger statement: if the initial topology is $T_0$, 
then the D-topology is completely Hausdorff, 
\emph{i.e.}, points are separated by continuous functions $X\to[0,1]$
(in fact, by diffeologically smooth functions contained in $\calF$).
\end{itemize}

}

\begin{appendices}

\section{Proofs}
\labell{proofs}

\subsection{Reflexive stability}
\labell{subsec:reflexive stability}

Recall that, given a set $X$ with a collection $\calD_0$
of parametrizations and a collection $\calF_0$ of real-valued functions,
$\Phi \calD_0$ denotes the set of those real-valued functions
$f \colon X \to \R$ whose precomposition with each element of $\calD_0$
is infinitely-differentiable,
and $\Pi \calF_0$ denotes the set of those parametrisations
$p \colon U \to X$
whose composition with each element of $\calF_0$ 
is infinitely-differentiable.

\begin{lemma}\labell{l:phidiffstr}
Fix a set $X$, and let $\mathcal{D}_0$ be a family of parametrisations
into $X$.
Then $\Phi\mathcal{D}_0$ is a differential structure on $X$.
\end{lemma}

\begin{proof}
We first show smooth compatibility.
Let $f_1,...,f_k\in\Phi\mathcal{D}_0$ and let $F\in\CIN(\RR^k)$.
Let $p \in \calD_0$.  
Because the components of $(f_1,\ldots,f_k) \circ p$ 
are infinitely-differentiable, 
the composition $F \circ (f_1,\ldots,f_k) \circ p$ is infinitely-differentiable.
Because $p$ is arbitrary, $F \circ (f_1,\ldots,f_k)$ is in $\Phi\calD_0$.

We now show locality.  
Equip $X$ with the initial topology of $\Phi\calD_0$.
Let $f \colon X\to\RR$ be a function satisfying:
for every $x \in X$ there is an open neighbourhood $V$ of $x$ in $X$
and a function $g\in\Phi\mathcal{D}_0$ such that $f|_V=g|_V$.
We want to show that $f\in\Phi\mathcal{D}_0$.
Fix $(p \colon U\to X)\in\mathcal{D}_0$.  Let $V\subseteq X$ be an open subset, 
and let $g\in\Phi\mathcal{D}_0$ be a function such that $f|_V=g|_V$.
Then $f\circ p|_{p^{-1}(V)}=g\circ p|_{p^{-1}(V)}$.  
The pre-image $p^{-1}(V)$ is open in $U$. (Indeed, $V$ is a union of pre-images $h^{-1}((a,b))$ of open intervals $(a,b)$ under functions $h$ in $\Phi\calD_0$, so $p^{-1}(V)$ is a union of the pre-images $(h\circ p)^{-1}((a,b))$, and $h\circ p\colon U\to\RR$ is infinitely-differentiable, hence continuous, because $p\in\calD_0$ and $h\in\Phi\calD_0$.)  Since each such $g\circ p$ is smooth in $U$ and is covered by such open sets $p^{-1}(V)$, and since smoothness is a local condition,
$f\circ p \colon U \to \RR$ is smooth. 
Since $p \in \calD_0$ is arbitrary, $f\in\Phi\mathcal{D}_0$.
\end{proof}

\begin{lemma}\labell{l:pidiffeol}
Fix a set $X$, and let $\calF_0$ be a set of real-valued functions on $X$.
Then $\Pi\mathcal{F}_0$ is a diffeology on $X$.
\end{lemma}

\begin{proof}
To see that $\Pi\mathcal{F}_0$ contains all the constant maps into $X$,
note that if $p \colon U \to X$ is constant 
then for any $f\in\mathcal{F}_0$ the composition $f\circ p \colon U \to X$ 
is constant, hence infinitely-differentiable.

Next, we show locality. Let $p \colon U\to X$ be a parametrisation
such that for every $u\in U$ there is an open neighbourhood $V$
of $u$ in $U$ such that $p|_V\in\Pi\mathcal{F}_0$;
we want to show that $p\in\Pi\mathcal{F}_0$.
Let $f\in\mathcal{F}_0$.
For any $u\in U$, there is an open neighbourhood $V$ of $u$ in $U$
such that $f\circ p|_V$ is smooth.
Since smoothness on $U$ is a local condition, 
$f\circ p \colon U \to \RR$ is smooth.
Since $f \in \calF_0$ is arbitrary, $p\in\Pi\mathcal{F}_0$.

Finally, we show smooth compatibility.  Let $U$ and $V$ be open
subsets of cartesian spaces, and let $F \colon V\to U$ be a smooth map.
Let $(p \colon U\to X)\in\Pi\mathcal{F}_0$.
For any $f\in\mathcal{F}_0$, we have that $f \circ p$ is smooth,
so $f\circ p\circ F$ is smooth.
Because $f \in \calF_0$ is arbitrary, $p\circ F\in\Pi\mathcal{F}_0$.
\end{proof}

\begin{proof}[Proof of Reflexive Stability (Proposition~\ref{p:reflexive})]
By Lemma~\ref{l:phidiffstr}, $\calF := \Phi\calD_0$ is a differential structure;
by Remark \ref{properties}, it is reflexive.
By Lemma~\ref{l:pidiffeol}, $\calD := \Pi \calF_0$ is a diffeology;
by Remark \ref{properties}, it is reflexive.
\end{proof}

\subsection{Isomorphism of categories of reflexive spaces} 
\labell{subsec:proof of A}

Recall that $\bfPhi(X,\calD) = (X,\Phi\calD)$ on objects
and $\bfPhi(F)=F$ on morphisms.

\smallskip\noindent\emph{Proof that $\bfPhi$ is a functor 
from the category of diffeological spaces
to the category of reflexive differential spaces.}
By Proposition~\ref{p:reflexive},
if $(X,\mathcal{D})$ is a diffeological space
then $(X,\Phi\mathcal{D})$ is a reflexive differential space.
We need to show that if
$F \colon (X,\mathcal{D}_X)\to(Y,\mathcal{D}_Y)$
is diffeologically smooth 
then $F$ is also functionally smooth as a map between
the reflexive differential spaces
$(X,\Phi\mathcal{D}_X)$ and $(Y,\Phi\mathcal{D}_Y)$.
Let $f\in\Phi\mathcal{D}_Y$.
Let $p \in \calD_X$.  Because $F$ is diffeologically smooth,
$F \circ p \in \calD_Y$.
This and the fact that $f \in \Phi\calD_Y$
imply that $f \circ F \circ p$ is infinitely-differentiable.
Since $p \in \calD_X$ is arbitrary,
this shows that $f \circ F \in \Phi\calD_X$.
Since $f \in \Phi\calD_Y$ is arbitrary,
this shows that $F$ is functionally smooth.
\qed

Recall that $\bfPi(X,\calF) = (X,\Pi\calF)$ on objects
and $\bfPi(F)=F$ on morphisms.

\smallskip\noindent\emph{Proof that $\bfPi$ is a functor 
from differential spaces to reflexive diffeological spaces.}
By Proposition~\ref{p:reflexive},
if $(X,\mathcal{F})$ is a differential space,
then $(X,\Pi\mathcal{F})$ is a reflexive diffeological space.
We need to show that if
$F \colon (X,\mathcal{F}_X)\to(Y,\mathcal{F}_Y)$
is functionally smooth 
then $F$ is also diffeologically smooth as a map between
the reflexive diffeological spaces $(X,\Pi\mathcal{F}_X)$ and
$(Y,\Pi\mathcal{F}_Y)$.  Let $p\in\Pi\mathcal{F}_X$.
Let $f \in \calF_Y$.  Because $F$ is functionally smooth,
$f \circ F \in \calF_X$.  This and the fact that $p \in \Pi\calF_X$
imply that $f \circ F \circ p$ is smooth.
Since $f \in \calF_Y$ is arbitrary, this shows
that $F \circ p \in \Pi\calF_Y$.
Because $p \in \Pi\calF_X$ is arbitrary, this shows that $F$
is diffeologically smooth.
\qed

\begin{proof}[Proof of isomorphism of categories of reflexive spaces
(Theorem \ref{t:A})]
If $(X,\mathcal{F})$ is a reflexive differential space, then
$ \bfPhi \circ \bfPi (X,\calF) = (X,\Phi\Pi\calF) = (X,\calF)$.
If $(X,\mathcal{D})$ is a reflexive diffeological space,
then $\bfPi \circ \bfPhi (X,\calD) = (X,\Pi\Phi\calD) = (X,\calD)$.
This and the fact that $\bfPi$ and $\bfPhi$ send every map to itself
shows that the restriction of the functor $\mathbf{\Phi}$
to the subcategory of reflexive diffeological spaces
and the restriction of the functor $\mathbf{\Pi}$
to the subcategory of reflexive differential spaces
are inverses of each other and give an isomorphism of categories.
\end{proof}

\subsection{Fr\"olicher spaces as reflexive spaces}
\labell{subsec:proof of B}

Recall that $\bfXi(X,\calC,\calF) = (X,\calF)$ on objects
and $\bfXi(F) = F$ on morphisms.

\smallskip

\noindent\emph{Proof that $\bfXi$ is a functor from the category
of Fr\"olicher spaces to the category of reflexive diffeological spaces.}
Let $(X,\calC,\calF)$ be a Fr\"olicher space.
In particular, $\calF = \Phi\calC$.
By Proposition~\ref{p:reflexive}, 
$\calF$ is a reflexive differential structure.  
Thus, $\bfXi$ takes Fr\"olicher spaces to reflexive differential spaces.
As noted in Definition~\ref{d:frolicher},
if a map of Fr\"olicher spaces is Fr\"olicher smooth,
then it is also functionally smooth.
\qed

\smallskip

Recall that
$ \bfGamma(X,\calF) = (X,\Gamma\calF,\Phi\Gamma\calF) $
on objects and $\bfGamma(F) = F$ on morphisms.

\smallskip\noindent\emph{Proof that $\bfGamma$ is a functor 
from the category of differential spaces
to the category of Fr\"olicher spaces.}
Let $(X,\calF)$ be a differential space.
The equality $\Gamma \Phi \Gamma \calF = \Gamma \calF$
shows that $(X,\Gamma\calF,\Phi\Gamma\calF)$ is a Fr\"olicher space.
As noted in Definition~\ref{d:frolicher},
if a map of differential spaces is functionally smooth,
then it is also Fr\"olicher smooth.
\qed

\begin{proof}[Proof of ``Fr\"olicher spaces as reflexive spaces''
(Theorem \ref{t:B})]
If $f$ is a real-valued function on a diffeological space $(X,\calD)$ and $f\circ c$ is infinitely-differentiable for every plot $c$ in $\calD$ with domain $\RR$, then $f\circ p$ is infinitely-differentiable for every plot $p\colon U\to X$ in $\calD$.  Indeed, 
by Boman's theorem \cite[Theorem 1]{boman}
it is enough to show that the composition $f \circ p \circ \gamma$
is infinitely-differentiable for every infinitely-differentiable curve
$\gamma \colon \RR \to U$, and this is true because $p \circ \gamma$ is a plot in $\calD$
with domain $\RR$. 

If $(X,\calC,\calF)$ is a Fr\"olicher space, 
then $\bfGamma \circ \bfXi (X,\calC,\calF) = \bfGamma (X,\calF)
 = (X,\Gamma\calF,\Phi\Gamma\calF) = (X,\calC,\calF)$.
If $(X,\calF)$ is a reflexive differential space, 
then $\bfXi \circ \bfGamma(X,\calF) = \bfXi(X,\Gamma\calF,\Phi\Gamma\calF)
 = (X,\Phi\Gamma\calF) = (X,\Phi\Pi\calF) = (X,\calF)$.
Here, the equality $\Phi\Gamma\calF = \Phi\Pi\calF$ is obtained
from the previous paragraph by setting $\calD = \Pi\calF$.
This and the fact that $\bfPi$ and $\bfGamma$ send every map to itself
shows that the functor $\bfXi$ and the restriction of the functor $\bfPi$
to the category of reflexive differential spaces 
are inverses of each other and give an isomorphism of categories.
\end{proof}

\subsection{Intersecting submanifolds}
\labell{subsec:lines in plane}
%

\begin{proof}[Proof of Part (1) of Example~\ref{x:wedgesum}]

Recall that $E \subset \R^2$ is the union of the two coordinate axes
and $\CIN(E)$ is its subspace differential structure,
that $\calD_X$ is the quotient diffeology on 
$X := ( \RR_1 \amalg \RR_2) / (0_1 \sim 0_2)$,
and that $\varphi \colon X \to E$ 
is the bijection whose pullback to $\R_1$ is $x \mapsto (x,0)$
and whose pullback to $\R_2$ is $y \mapsto (0,y)$.
Fix a real-valued function $f \colon E \to \RR$.
Define $f_i \colon \RR \to \RR$, for $i=1,2$, by
$f_1(x) = f(x,0)$ and $f_2(y) = f(0,y)$.
We need to show that each of the conditions $f \in \CIN(E)$
and $\varphi^* f \in \Phi \calD_X$
is equivalent to $f_1$ and $f_2$ being smooth.

First, suppose that $f \in \CIN(E)$.
Then $f_1$ and $f_2$, being the compositions of the smooth maps 
$x \mapsto (x,0)$ and $y \mapsto (0,y)$ with a smooth extension
of $f$ to $\RR^2$, are smooth.

Now, suppose that $\varphi^* f \in \Phi \calD_X$.
Let $i=1$ or $i=2$.
The inclusion map of the $i$th copy of $\RR$ in $X$,
which we denote $I_i \colon \RR \to X$,
is in the quotient diffeology $\calD_X$.
By the definition of $\Phi \calD_X$,
the composition $(\varphi^* f) \circ I_i$ is smooth.
This composition is $f_i$, so $f_i$ is smooth.

Now, suppose that $f_1$ and $f_2$ are smooth.
Then $(x,y) \mapsto f_1(x) + f_2(y) - f(0,0)$
is a smooth extension of $f$ to $\RR^2$.
This shows that $f \in \CIN(E)$.

Still assuming that $f_1$ and $f_2$ are smooth,
let $p \colon U \to X$ be a plot in the quotient diffeology~$\calD_X$.
Let $u \in U$ be any point.
Let $V$ be a connected neighbourhood of $u$ in $U$
and $\tilde{p} \colon V \to \RR_1 \amalg \RR_2$ a smooth lifting 
of $p|_V$;
these exist by the definition of the quotient diffeology.
By continuity, the image of $\tilde{p}$ is contained in $\RR_i$
for some $i \in \{1,2\}$.
The map $(\varphi^* f) \circ p|_V \colon V \to \RR$,
being the composition of the smooth maps $\tilde{p}$ and $f_i$, is smooth.
Since smoothness is a local condition and $u \in U$ is arbitrary,
$(\varphi^* f) \circ p \colon U \to \RR$ is smooth.
Since $p \in \calD_X$ is arbitrary, $\varphi^* f \in \Phi\calD_X$.
\end{proof}

\begin{proof}[Sketch of proof of Example~\ref{x:wedgesum2}]
The proof is similar to that of Example~\ref{x:wedgesum} once we make 
the following observation.
Let $\calI$ be a set of subsets of $\{ 1, \ldots, n \}$,
and let $E_\calI := \bigcup\limits_{I \in \calI} \R^I$ where $\RR^I$ is the span of the $x_i$-axes for $i\in I$.
For every $I \in \calI$, 
let $\pr_I \colon \RR^n \to \RR^I$ denote the natural projection map.
Then, for every function $f \colon E_\calI \to \RR$
that is smooth on $\RR^I$ for each $I \in \calI$, the function
$$ \sum\limits_{\substack{A \subset \calI \\ A \neq \emptyset}}
 (-1)^{1+|A|} f \circ \pr_{\cap A} \colon \RR^n \to \RR $$
is a a smooth extension of $f$ to $\R^n$.
Indeed, to see that this function coincides with $f$ on $E_\calI$, 
we argue as follows.  
Fix $x \in E_\calI$.  When $A$ is the empty set, then $\pr_{\cap A}(x) = x$
(this is not the same as when $\bigcap A=\emptyset$, in which case $\pr_\emptyset=0$).
Hence $(-1)^{1+|A|} f \circ \pr_{\cap A}(x) = -f(x)$. 
So it is enough to show that the sum
$$ \sum\limits_{A \subset \calI} (-1)^{1+|A|} f(\pr_{\cap A}(x)) $$
vanishes.  We can write this sum as
$$ \sum\limits_{\substack{A \subseteq \{I_1,\ldots,I_m\} \\ 
 B\subseteq \{ J_1,\ldots,J_s \}}} (-1)^{1+|A|+|B|} 
f (\pr_{(\cap A) \cap (\cap B)}(x)) ,$$
where $I_1, \ldots, I_m$ is an enumeration of the set
$\{ I \in \calI \ | \ x \in \R^I \}$ 
and $J_1,\ldots,J_s$ is an enumeration of the set
$\calI \setminus \{ I_1,\ldots,I_m \}$.
Because $x \in \R^{I_i}$ for all $i=1,\ldots,m$, we have
$\pr_{ (\cap A) \cap (\cap B) } (x) = \pr_{\cap B}(x)$. 
So we can write the above sum as
$$ \sum\limits_{ B \subseteq \{ J_1 , \ldots , J_s \} } 
 (-1)^{1+|B|}  f(\pr_{\cap B}(x))
 \sum\limits_{ A \subseteq \{ I_1,\ldots, I_m\} } (-1)^{|A|} .$$
Because $x \in E_\calI$, we have $m \geq 1$,
and so $\sum\limits_{A \subseteq \{ I_1,\ldots,I_m \} } (-1)^{|A|} = 0$,
so the above sum vanishes.
\end{proof}

\begin{proof}[Completion of the proof 
of Part (1) of Example~\ref{three lines}]

Recall that $E \subseteq \RR^3$ is the union of the three coordinate axes;
$S \subseteq \RR^2$ is the union $l_1 \cup l_2 \cup l_3$
where $l_1$ is the $x$-axis, $l_2$ is the $y$-axis,
and $l_3$ is the line given by $y = x$;
and $\varphi \colon E \to S$ is the map
$(t,0,0) \mapsto (t,0)$, $(0,t,0) \mapsto (0,t)$, 
$(0,0,t) \mapsto (t, t)$. 
We need to prove that $\calD_S \subseteq \varphi \circ \calD_E$.
For this, we fix an open subset $U$ of $\RR^k$ for some $k$
and a plot 
$$ p \colon U \to S $$ 
of $S$, and we need to prove that $\varphi^{-1} \circ p \colon U \to E$
is a plot of $E$.
Let $(p_1,p_2) \colon U \to \RR^2$ be the composition of $p \colon U \to S$
with the inclusion map $S \to \RR^2$,
and let $q \colon U \to \R^3$ be the composition 
of $\varphi^{-1} \circ p \colon U \to E$
with the inclusion map $E \to \RR^3$.  
On each subset $p^{-1}(l_i)$, the map $q$ 
coincides with the map $g_i$, where
$$ g_1(u) = (p_1(u),0,0) , \quad
   g_2(u) = (0,p_2(u),0) , \quad \text{ and } \quad
   g_3(u) = (0,0,p_1(u)) . $$
The maps $g_i \colon U \to \RR$ are smooth (because $p$ is a plot),
and we need to prove that the map $q \colon U \to \R^3$ is smooth.

Let
$$ U_i = \text{interior}( p^{-1}(l_i) ) 
\ \ \text{ for } i=1,2,3, \quad \text{ and let } \quad
   W = \bigcup_{j \neq k} \ol{U_j} \cap \ol{U_k} .  $$
We claim that
\begin{equation} \labell{union W}
U = U_1 \cup U_2 \cup U_3 \cup W ,
\quad \text{ and } \quad 
\ol{U_i} \subseteq U_i \cup W \ \ \text{ for } i=1,2,3.
\end{equation}
Indeed, let $u \in U \ssminus (U_1 \cup U_2 \cup U_3)$.
Then $p(u) = 0$ and each neighbourhood of $u$
contains points from at least two of the sets
$p^{-1}(l_i \ssminus \{ 0 \} )$ for $i = 1,2,3$.
So there exist $j \neq k$ such that every neighbourhood of $u$ 
contains points of $p^{-1}(l_j \ssminus \{ 0 \} )$ 
and points of $p^{-1}(l_k \ssminus \{ 0 \} )$.
Then $u \in \ol{U_j} \cap \ol{U_k}$, and so $u \in W$.
This proves the first part of~\eqref{union W}.
Now suppose that $u \in \ol{U_i}$.
By the first part of~\eqref{union W}, 
either $u \in U_i$, or $u \in W$, or $u \in U_j$ for $j \neq i$.
In the first or second case, $u \in U_i \cup W$.
In the third case, $u \in \ol{U_i} \cap U_j 
\subseteq \ol{U_i} \cap \ol{U_j} \subseteq W \subseteq U_i \cup W$.
This proves the second part of~\eqref{union W}.

Let $t_1,\ldots,t_k$ be the coordinates on $U \subseteq \RR^k$.
Consider the differentiation operators
$$ D_m = \frac{ \del^{m_1+\ldots+m_k} }
        { \del t_1^{m_1} \cdots \del t_k^{m_k} } 
\quad \text{ for }  m = (m_1,\ldots,m_k) \in \ZZ_{\geq 0}^k .$$
For each $i \in \{ 1,2,3 \}$ the restriction $D_m (p_1,p_2)|_{U_i}$
takes values in the linear subspace $l_i$ of $\R^2$.
By continuity, $(D_m (p_1,p_2))|_{\ol{U_i}}$ also takes values in $l_i$.
If $j \neq k$, then,  because $l_j \cap l_k = \{ 0 \}$,
the derivatives $D_m(p_1,p_2)$ vanish on $\ol{U_j} \cap \ol{U_k}$.  So
\begin{equation} \labell{derivatives vanish}
\text{if $u \in W$, then } 
 D_m g_i (u) = 0 \ 
   \text{ for all } m \in \Z_{\geq 0}^k \text{ and } i=1,2,3.
\end{equation}
Consider the following statements.
\begin{itemize}[noitemsep,topsep=0pt]
\item[(I$_m$)] $D_m q \colon U \to \RR^3$ exists throughout $U$
and vanishes on $W$.
\item[(II$_m$)]
For each $i=1,2,3$, $D_m q$ (exists and) coincides with $D_m g_i$ on $\ol{U_i}$.
\item[(III$_m$)]
$D_m q \colon U \to \RR^3$ (exists and) is continuous.
\end{itemize}

(I$_m$) implies (II$_m$).
This follows by the second part of~\eqref{union W}
from the facts that $q$ and $g_i$ coincide
on the open set $U_i$ and that, assuming (I$_m$),
 $D_m q$ and $D_mg_i$ both vanish
at the points of $W$ ($D_m q$ by hypothesis and $D_mg_i$ 
by \eqref{derivatives vanish}).

(I$_m$) and (II$_m$) imply (III$_m$).
This is because $D_m q$ coincides with continuous maps
on the closed sets $\ol{U_1}$, $\ol{U_2}$, $\ol{U_3}$, $W$,
whose union is $U$ (by the first part of~\eqref{union W}).

We will now show that (I$_m$) is true for all $m$.
For $m=0$, this follows from \eqref{derivatives vanish}.
Arguing by induction, assume that (I$_{m'}$), 
and hence (II$_{m'}$) and (III$_{m'}$), are true,
and let $m$ be obtained from $m'$
by increasing one of its coordinates by one, say, the $\ell$th coordinate.
Because $q$ coincides with the smooth map $g_i$ on the open set $U_i$,
the derivative $D_mq$ exists on the $U_i$s.
Denote the $\ell$th standard basis element of $\R^k$ by $e_\ell$. 
Fix a point $u \in W$.  
For any $h$ such that $u+he_\ell \in U$, we claim that 
\begin{equation} \labell{difference quotient}
 \frac{ D_{m'}q (u+he_\ell) - D_{m'}q(u) }{ h }
 = \begin{cases}
\displaystyle{
         \frac{ D_{m'}g_i (u+he_\ell) - D_{m'}g_i (u) }{ h } } 
                        & \quad \text{ if } u+he_\ell \in U_i \\[10pt]
         \ \hfill 0 \hfill \  & \quad \text{ if } u+he_\ell 
                  \not\in U_1 \cup U_2 \cup U_3 .
\end{cases}
\end{equation}
The first case is because 
$q$ and $g_i$, and hence their derivatives, coincide 
on the open subset $U_i$,
and because $D_{m'}q(u) = 0$ (by (I$_{m'}$))
and $D_{m'}g_i(u) = 0$ (by~\eqref{derivatives vanish} for $m'$).
In the second case $u+he_\ell \in W$
(by the first part of \eqref{union W})
and $u \in W$ (by assumption),
so $D_{m'}q (u+he_\ell) = D_{m'}q (u) = 0$ (by (I$_{m'}$)).
Since each term on the right hand side of~\eqref{difference quotient}
converges to zero as $h \to 0$ (by \eqref{derivatives vanish} for $m$),
we conclude that the left hand side converges to zero,
so $D_m q(u)$ exists and is equal to zero.
Because $u \in W$ is arbitrary, we obtain~(I$_m$).

Thus, (I$_m$), (II$_m$), and (III$_m$) are true for all $m$.
In particular, $q$ is smooth, as required.
\end{proof}

\section{Comparisons with Other Structures}\labell{app:comparisons}

In this appendix, we compare diffeological and differential spaces 
(and hence Fr\"olicher spaces by Theorem~\ref{t:B}) with some of the other 
generalisations of differential calculus that appear in the literature.
We refer to Andrew Stacey's paper \cite{stacey}
for a more extensive comparative study of Chen spaces, Smith spaces,
diffeological spaces, Fr\"olicher spaces, and differential spaces; we do not address Chen spaces nor Smith spaces here.  For a direct comparison of diffeological and Chen spaces, see \cite{KW:convex}.  We refer to Joao Nuno Mestre's Ph.D.\ thesis \cite[Chapter 2]{mestre} for a comparative study of differential spaces, Mostow spaces, subcartesian spaces, differentiable spaces (not the same as ``differential'' spaces), and $\CIN$-schemes.

For the sake of brevity, we do not give definitions of the structures discussed below. Instead, we refer the reader to the following sources that are more focused on these subjects (we do not claim that this is an exhaustive list):
\begin{itemize}
	\item Lie groupoids and stacks \cite{BX,lerman,MM,pronk};
	\item sheaves of sets over a site \cite{BH};
	\item synthetic differential geometry and $\CIN$-schemes \cite{dubuc, joyce, kock, MR};
	\item Mostow spaces \cite{mostow};
	\item subcartesian spaces \cite{A,AS1,AS2,sniatycki};
	\item differentiable spaces \cite{GS,spallek1,spallek2,spallek3}.
\end{itemize}

\subsection{From Lie groupoids to diffeological spaces}
\labell{ss:groupoids}

There is a functor from the category of Lie groupoids to the category
of diffeological spaces, sending a Lie groupoid to its orbit space
equipped with the quotient diffeology, and smooth morphisms between
Lie groupoids to diffeologically smooth maps between the orbit spaces.
This functor is neither faithful nor full, even when restricted to
effective \'etale proper Lie groupoids (\emph{i.e.}\ effective orbifolds);
see the examples of Moshe Zadka \cite[Examples 24 and 25]{IZKZ}.

However, this functor does factor through Morita equivalence.  In fact,
the bicategory of Lie groupoids with bibundles between them
and isomorphisms of bibundles as 2-arrows (see \cite{lerman,MM} for definitions) has a pseudofunctor to
diffeological spaces, in which bibundles are sent to diffeologically
smooth maps, and 2-arrows are sent to trivial 2-arrows (diffeological
spaces form an honest $1$-category); this is proven by Watts in \cite[Theorem 3.8]{watts-gpds}.
Moreover, when this pseudofunctor is restricted to effective orbifolds with
``locally invertible'' bibundles between them, then this is an equivalence
of categories onto diffeological orbifolds with locally invertible smooth
maps between them.  In fact, Karshon and Miyamoto in \cite{quasifolds} 
(following an earlier preprint by Karshon and Zoghi
that was announced in \cite{zoghi})
prove that this restriction works in the more general setting of so-called effective quasifold groupoids and diffeological quasifolds.

An important feature of this pseudofunctor is that much of the isotropic information
is generally lost.  Indeed, consider $\U(n)$ and $\SO(2n)$ acting
on $\RR^{2n}$ by rotations.  The resulting diffeological quotients
are diffeologically diffeomorphic, but the groupoids are not even
Morita equivalent: the stabilisers at the origin are not isomorphic.
However, in certain circumstances these stabilisers can be recovered: it follows from \cite{IZKZ} that the restricted pseudofunctor from effective \'etale proper Lie groupoids to diffeological orbifolds is injective on objects.  It follows from \cite{CDGMW} that the same pseudofunctor restricted to action groupoids of faithful linear representations of the circle is also injective on objects.  See Example~\ref{x:orthogonal action} for more such examples.  It would be interesting to pin
down precisely which isotropic information is lost and which is retained,
even in the case of linear compact group actions.  There are also examples of Lie groupoids whose orbit spaces are diffeologically diffeomorphic, the isotropic information is the same, but the Lie groupoids are not Morita equivalent; see \cite[Section 7]{quasifolds}.

\subsection{From stacks to diffeological spaces}
\labell{ss:stacks}

There is an equivalence of bicategories between the bicategory of Lie groupoids and the $2$-category of differentiable stacks,
and thus by \ref{ss:groupoids} there is a pseudofunctor from differentiable
stacks to diffeological spaces (viewing diffeological spaces as a bicategory with trivial $2$-arrows).  The image of a differentiable stack via this pseudofunctor is the orbit space of a Lie groupoid representing the stack, but different choices of representative Lie groupoid yield only diffeomorphic orbit spaces.

In Watts-Wolbert \cite{WW}, the authors show that this pseudofunctor can be described more strictly as a $2$-functor: given a differentiable stack, this $2$-functor sends it to a diffeological space that only depends on the stack (although it is diffeomorphic to the orbit space of any representative Lie groupoid).  Moreover, this $2$-functor extends to all stacks over
the site of smooth manifolds, sending a stack to its underlying
diffeological ``coarse moduli space''.  Up to an application of the
comparison lemma of sheaves on sites, the $2$-functor factors through
the so-called concretization functor of Baez and Hoffnung \cite{BH},
which sends a sheaf of sets over the site of smooth manifolds to its
underlying diffeological space.  In fact, the $2$-functor is adjoint to
the inclusion functor from diffeological spaces into stacks, which again
factors through sheaves of sets over the site of manifolds; see Subsection~\ref{ss:diffeol sheaves} for more details.  In other
words, stacks form a language that unifies Lie groupoids up to Morita
equivalence, sheaves of sets over manifolds, and diffeological spaces.

Differential spaces do not fit into this setting at all.  For instance, $\RR$ with its standard diffeology is diffeologically diffeomorphic to the cusp given by $x^2=y^3$ in $\RR^2$; see Karshon-Miyamoto-Watts \cite{KMW}.  However, the standard differential structure on $\RR$ is not functionally diffeomorphic to the subspace differential structure on the cusp; indeed, the structural dimension of the cusp is $2$, whereas at each non-cuspoidal point, it is $1$.  It is also equal to $1$ at all points of $\RR$.  Since structural dimension is an invariant of such differential spaces (see \cite{LSW}), these two spaces are not functionally diffeomorphic.

\subsection{Diffeological spaces as sheaves over $\Open$}
\labell{ss:diffeol sheaves}

A diffeology on a set $X$ cannot be obtained from a ``structure sheaf''
on $X$ as a topological space;
for example, the irrational torus (Example~\ref{x:irrational flow}) 
has an interesting diffeology but a trivial topology.
Instead, diffeology can be viewed as a sheaf of sets over a site.

Namely, let $\Open$ denote the category whose objects
are the open subsets of cartesian spaces ($\R^n$, $n \geq 0$)
and whose arrows are smooth maps.  This is a site whose coverages are exactly the standard open covers of open subsets of cartesian spaces.
A diffeology $\calD_X$ on a set $X$ 
determines a contravariant functor $\calD_X\colon\Open\to\Set^\op$ (\emph{i.e.}\ a presheaf), sending $U$ to the set of plots with domain $U$. The locality axiom of diffeology guarantees that this presheaf is in fact a sheaf.  Furthermore, if we let $\ul{X}$ be the sheaf assigning to each object $U$ of $\Open$ \emph{all} functions $U\to X$, then $\calD_X$ is a subsheaf of $\ul{X}$ that satisfies $\calD_X(\RR^0)=\ul{X}(\RR^0)$.  This definition is due to Lerman \cite[Definition~A.13]{parallel-transport-stacks}.  A similar approach to diffeology as sheaves over a category already appears in the work of Iglesias-Zemmour in the appendix to his 1986 paper \cite{iglesias:1986}.

In the language of Baez and Hoffnung, diffeologies are exactly the sheaves over $\Open$ that are ``concrete'' \cite{BH}.  Not all sheaves over $\Open$ are concrete, however.  For instance, for $k>0$, consider the sheaf $\Omega^k(\cdot)$ assigning to the object $U$ in $\Open$ all differential $k$-forms $\Omega^k(U)$.  Then $\Omega^k(\RR^0)$ is trivial.  Note, however, that this does not correspond to the diffeological space $\RR^0$, as $\calD_{\RR^0}$ sends $U$ to the singleton consisting of the constant plot $U\to\RR^0$, whereas $\Omega^k(U)$ is not a singleton for general $U$.

\subsection{Differential structures and ringed spaces}
\labell{ss:sheaves}

Let $(X,\calF)$ be a differential space.  Then there is a naturally
induced reduced ringed space $(X,\widehat{\calF})$ 
where for each
open $U\subseteq X$ (with respect to the initial topology induced
by $\calF$), the ring $\calF(U)$ is the subspace differential
structure on $U$.  (A ``reduced'' ringed space is a ringed space whose sheaf is a sheaf of continuous real-valued functions.)
Moreover, any smooth map of differential spaces
$f\colon(X,\calF_X)\to(Y,\calF_Y)$ induces a morphism of reduced ringed
spaces $(f,f^\sharp)\colon(X,\widehat{\calF}_X)\to(Y,\widehat{\calF}_Y)$.
One can recover
$(X,\calF)$ by taking $\calF$ to be the ring of global sections.

The opposite operation does not work: starting with an
appropriate reduced ringed space, the sheafification of the global
sections as above does not necessarily return the original ringed space.
For example, consider the non-Hausdorff manifold ``the real line with two origins''; see Example~\ref{e:double origin}.  This is the quotient of $\RR\amalg\RR$
by the relation $x\sim y$ if $x$ and $y$ are non-zero and copies of the
same real number.  Equip the resulting quotient topological space $X$ with the sheaf sending an open set $U$ (in the quotient topology) to the subspace differential structure on $U$ induced by the quotient differential structure on $X$.  This sheaf cannot be obtained from a differential space, as the topology on $X$ is not initial with respect to the global sections of the sheaf.

\subsection{From differential spaces to $\CIN$-schemes}
\labell{ss:c-schemes}

A $\CIN$-ringed space is a topological space equipped with a sheaf of
$\CIN$-rings.  An affine $\CIN$-scheme is a locally $\CIN$-ringed space
isomorphic to the real spectrum of a $\CIN$-ring.  A $\CIN$-scheme is
a locally $\CIN$-ringed space $(X,\calO_X)$ such that $X$ admits an
open cover $\{U_\alpha\}$ in which $(U_\alpha,\calO_X|_{U_\alpha})$
is an affine $\CIN$-scheme.

Proposition 2.77 of \cite{mestre} states that any reduced affine
$\CIN$-scheme can be considered to be a differential space, 
and Corollary of 2.78 of
\cite{mestre} states that any reduced $\CIN$-scheme can be considered to be a Mostow space.
In general, however, a reduced $\CIN$-scheme is not the sheafification
of a differential space; see the example in \ref{ss:sheaves}.  On the other hand, differential spaces embed into (reduced) $\CIN$-schemes; see \cite[Corollary 2.76]{mestre}.

We must keep the adjective ``reduced'' above.  Indeed, even for affine $\CIN$-schemes, there may be, for instance, nilpotent elements in the $\CIN$-ring of global sections that cannot be realised as real-valued functions on a set.  A simple example of this is given by the so-called ``dual numbers'': $\RR[x]/(x^2)$ has a real spectrum given by a point with corresponding sheaf of functions exactly those of $\RR^0$.  All elements generated by $x(x^2)$ are forgotten by this sheaf.  See \cite[page 19]{MR} for more details on this example.

\subsection{From differential spaces to Mostow spaces}
\labell{ss:mostow}

A Mostow space is a reduced ringed space $(X,\calF)$ such that for any open set
$U\subseteq X$, the ring $\calF(U)$ satisfies the smooth compatibility
condition of differential spaces: for any $f_1,\dots,f_k\in\calF(U)$ and any $g\in\CIN(\RR^k)$, the composition $g(f_1,\dots,f_k)$ is in $\calF(U)$.  It follows that there is a full
embedding of differential spaces into Mostow spaces by converting a
differential space into a ringed space as above.  Thus, one may view a
Mostow space as a differential space in which the topology on the space
is allowed to be finer than that of the initial topology.  For example,
given a diffeological space $(X,\calD)$ where $X$ is equipped with the
D-topology, the corresponding ringed space induced from $(X,\Phi\calD)$
is naturally a Mostow space.  This becomes a differential space if we
replace the D-topology with the initial topology.  In fact, replacing the topology of the Mostow space with the initial topology induced by its global sections is adjoint to the embedding of differential spaces into Mostow
spaces; see \cite[Proposition 2.66]{mestre}.

\subsection{From subcartesian spaces to differentiable spaces}
\labell{ss:subcartesian}

Recall from the introduction that a subcartesian space is a 
differential space that is locally functionally diffeomorphic to subsets
of cartesian spaces.  

Differentiable spaces are a special case (the ``$\infty$-standard'' case)
of spaces introduced by Spallek \cite{spallek1,spallek2,spallek3};
a standard reference on these is the book by Gonzalez-Salas
\cite{GS}.  A differentiable algebra is an $\RR$-algebra isomorphic
to $\CIN(\RR^n)/\mathfrak{a}$ for some $n$ and ideal $\mathfrak{a}$
closed with respect to the Fr\'echet topology.  An affine differentiable
space is a locally ringed space isomorphic to the real spectrum of
a differentiable algebra.  A differentiable space $(X,\calO_X)$ is
a locally ringed space admitting an open cover $\{U_\alpha\}$ such
that $(U_\alpha,\calO_X|_{U_\alpha})$ is an affine differentiable
space for each $\alpha$.  Warning: again, a \emph{differential} space and a
\emph{differentiable} space are two different things; we will prepend
``Sikorski'' to the former in this subsection to avoid confusion.

Proposition 2.81 of \cite{mestre} states that any reduced affine
differentiable space is a subcartesian space, and consequently,
any reduced differentiable space is a Mostow space.  These facts follow
from the discussion on $\CIN$-schemes mentioned above, and the fact
that differentiable spaces form a full subcategory of $\CIN$-schemes
\cite[Corollary 2.76]{mestre}.  In fact, one can say more.  It follows
from \cite[Proposition 5.6, Corollary 5.7]{GS} that any reduced affine
differentiable space is a closed Sikorski differential subspace of $\RR^n$
for some $n$, which is stronger than the subcartesian condition.

\subsection{Conclusion}
\labell{ss:conclusion}

It is the opinion of the author Watts that these
various theories described above should not be viewed as competing with each other,
but instead, that each theory individually focuses on specific attributes desired in
what one calls a ``smooth space''.  For example, neither diffeological,
Fr\"olicher, nor differential spaces start with a topological space; they
all start with structures consisting of functions mapping into and/or
out of a set. (While a differential structure uses the
initial topology to define the locality axiom, this topology is induced by the differential structure).  On the contrary, ringed spaces
such as $\CIN$-schemes, differentiable spaces, and Mostow spaces, start
with a topological space upon which a sheaf is defined.

As another example, Lie groupoids and stacks may encode more information than diffeological spaces.  In some situations, one can get away with using the simpler language of diffeology, however, in other situations, one needs the language of Lie groupoids and stacks.

A more unifying setting that keeps to concrete $1$-categories is to consider sets equipped with a diffeology \emph{and} a differential structure that are compatible; see Definition~\ref{d:compatible induce}.  This setting is used in \cite{virgin} and \cite{KW:convex}, and has the advantage of containing as full subcategories the categories of diffeological spaces, differential spaces, and Fr\"olicher spaces.  It also allows one to simultaneously utilise invariants and techniques specifically designed for diffeological or differential spaces.

\end{appendices}

{\ifworkmode
\newpage
\else
\end{document}